\documentclass[11pt,a4paper,reqno]{amsart} 
\usepackage{amsfonts,amsthm, amscd, epsfig, amsmath, amssymb,enumerate}
\def \beq {\begin{eqnarray}}
\def \eeq {\end{eqnarray}}
\def \beqn {\begin{eqnarray*}}
\def \eeqn {\end{eqnarray*}}

\newcommand{\halmos}{\rule{1ex}{1.4ex}}
\newcounter{for}[section]

\newtheorem{itlemma}{Lemma}[section]
\newtheorem{itproposition}[itlemma]{Proposition}
\newtheorem{itfact}[itlemma]{Fact}
\newtheorem{theorem}[itlemma]{Theorem}
\newtheorem{itcorollary}[itlemma]{Corollary}
\newtheorem{itremark}[itlemma]{Remark}
\newtheorem{itremarks}[itlemma]{Remarks}
\newtheorem{itdefinition}[itlemma]{Definition}
\newtheorem{itexample}[itlemma]{Example}

\newenvironment{fact}{\begin{itfact}\rm}{\end{itfact}}
\newenvironment{claim}{\begin{itclaim}\rm}{\end{itclaim}}
\newenvironment{lemma}{\begin{itlemma}}{\end{itlemma}}
\newenvironment{remark}{\begin{itremark}\rm}{\end{itremark}}
\newenvironment{remarks}{\begin{itremarks} \rm}{\end{itremarks}}
\newenvironment{corollary}{\begin{itcorollary}}{\end{itcorollary}}
\newenvironment{proposition}{\begin{itproposition}}{\end{itproposition}}
\newenvironment{definition}{\begin{itdefinition}\rm}{\end{itdefinition}}
\newenvironment{example}{\begin{itexample}\rm}{\end{itexample}}
\newcommand{\be}[1]{\addtocounter{for}{1} \begin{equation}\label{#1}}
\newcommand{\ee}{\end{equation}}
\newcommand{\bl}[1]{\begin{lemma}\label{#1}}
\newcommand{\br}[1]{\begin{remark}\label{#1}}
\newcommand{\brs}[1]{\begin{remarks}\label{#1}}
\newcommand{\bt}[1]{\begin{theorem}\label{#1}}
\newcommand{\bd}[1]{\begin{definition}\label{#1}}
\newcommand{\bp}[1]{\begin{proposition}\label{#1}}
\newcommand{\bfact}[1]{\begin{fact}\label{#1}}
\newcommand{\bc}[1]{\begin{corollary}\label{#1}}
\newcommand{\bex}[1]{\begin{example}\label{#1}}
\newcommand{\ec}{\end{corollary}}
\newcommand{\efact}{\end{fact}}
\newcommand{\eex}{\end{example}}
\newcommand{\el}{\end{lemma}}
\newcommand{\er}{\end{remark}}
\newcommand{\ers}{\end{remarks}}
\newcommand{\et}{\end{theorem}}
\newcommand{\ed}{\end{definition}}
\newcommand{\ep}{\end{proposition}}
\newcommand{\epr}{\end{proof}}
\newcommand{\bpr}{\begin{proof}}
\newcommand{\bcl}[1]{\begin{claim}\label{#1}}
\newcommand{\ecl}{\end{claim}}

\newcommand{\ecs}{\end{corollary}}
\newcommand{\eers}{\end{exercise}}
\newcommand{\eexs}{\end{example}}
\newcommand{\eems}{\end{example}}
\newcommand{\els}{\end{lemma}}
\newcommand{\eles}{\end{lemmaex}}
\newcommand{\ets}{\end{theorem}}
\newcommand{\eds}{\end{definition}}
\newcommand{\eps}{\end{proposition}}
\newcommand{\bi}{\begin{itemize}}
\newcommand{\ei}{\end{itemize}}
\newcommand{\ben}{\begin{enumerate}}
\newcommand{\een}{\end{enumerate}}

\def\vbar{\mathchoice{\vrule height6.3ptdepth-.5ptwidth.8pt\kern-.8pt}
   {\vrule height6.3ptdepth-.5ptwidth.8pt\kern-.8pt}
   {\vrule height4.1ptdepth-.35ptwidth.6pt\kern-.6pt}
   {\vrule height3.1ptdepth-.25ptwidth.5pt\kern-.5pt}}
\def\fudge{\mathchoice{}{}{\mkern.5mu}{\mkern.8mu}}
\def\bbc#1#2{{\rm \mkern#2mu\vbar\mkern-#2mu#1}}
\def\bbb#1{{\rm I\mkern-3.5mu #1}}
\def\bba#1#2{{\rm #1\mkern-#2mu\fudge #1}}
\def\bb#1{{\count4=`#1 \advance\count4by-64 \ifcase\count4\or\bba A{11.5}\or
   \bbb B\or\bbc C{5}\or\bbb D\or\bbb E\or\bbb F \or\bbc G{5}\or\bbb H\or
   \bbb I\or\bbc J{3}\or\bbb K\or\bbb L \or\bbb M\or\bbb N\or\bbc O{5} \or
   \bbb P\or\bbc Q{5}\or\bbb R\or\bbc S{4.2}\or\bba T{10.5}\or\bbc U{5}\or
   \bba V{12}\or\bba W{16.5}\or\bba X{11}\or\bba Y{11.7}\or\bba Z{7.5}\fi}}

\def \qed {{\hspace*{\fill}$\halmos$\medskip}}

\def \R {{\mathbb R}}

\def \N {{\mathbb N}}

\def \ra {\rightarrow }

\def \a{\alpha}

\def\g{\gamma}
\def\d{\delta}
\def\e{\varepsilon}

\def\b{\beta}

\def\l{\lambda}

\def\1{{\bf 1}}

\numberwithin{equation}{section}

\newtheorem{teo}{Theorem}[section]

\newtheorem{de}[teo]{Definition}

\renewcommand{\tilde}{\widetilde}

\def\proof{\noindent{\bf Proof.}\ }

\DeclareMathSymbol{\leqslant}{\mathalpha}{AMSa}{"36} % nicer `smaller or equal' 
\DeclareMathSymbol{\geqslant}{\mathalpha}{AMSa}{"3E} % nicer `larger or equal' 
\DeclareMathSymbol{\eset}{\mathalpha}{AMSb}{"3F}     % nicer `emptyset' 
\renewcommand{\leq}{\;\leqslant\;}                   % redef. of < or = 
\renewcommand{\geq}{\;\geqslant\;}                   % redef. of > or = 
\newcommand{\grad}{\nabla} 

\newcommand{\la}{\label} 
\def\1{\ifmmode {1\hskip -3pt \rm{I}} \else {\hbox {$1\hskip -3pt \rm{I}$}}\fi}

\newcommand{\cD}{\ensuremath{\mathcal D}} 
\newcommand{\cE}{\ensuremath{\mathcal E}}

\newcommand{\cL}{\ensuremath{\mathcal L}}

\newcommand{\bbN}{{\ensuremath{\mathbb N}} }

\newcommand{\bbR}{{\ensuremath{\mathbb R}} }

\newcommand{\bbZ}{{\ensuremath{\mathbb Z}} }

\newcommand{\ent}{{\rm Ent} } 
\newcommand{\var}{{\rm Var} }

\let\a=\alpha \let\b=\beta   \let\d=\delta  \let\e=\varepsilon
 \let\g=\gamma     \let\k=\kappa  \let\l=\lambda

   \let\G=\Gamma

\def\\{\hfill\break}

\def\thsp{\thinspace}

\def\tthsp{\kern .083333 em}

\def\?{\mskip -10mu}
\def\indbox#1{\hbox to \parindent{\hfil\ #1\hfil} }

\newfam\msafam
\newfam\msbfam
\newfam\eufmfam
\def\hexnumber#1{%
  \ifcase#1 0\or 1\or 2\or 3\or 4\or 5\or 6\or 7\or 8\or
  9\or A\or B\or C\or D\or E\or F\fi}
\font\tenmsa=msam10
\font\sevenmsa=msam7
\font\fivemsa=msam5
\textfont\msafam=\tenmsa
\scriptfont\msafam=\sevenmsa
\scriptscriptfont\msafam=\fivemsa        
\edef\msafamhexnumber{\hexnumber\msafam}%
\mathchardef\restriction"1\msafamhexnumber16
\mathchardef\ssim"0218
\mathchardef\square"0\msafamhexnumber03
\mathchardef\eqd"3\msafamhexnumber2C
\def\QED{\ifhmode\unskip\nobreak\fi\quad
  \ifmmode\square\else$\square$\fi}            
\font\tenmsb=msbm10
\font\sevenmsb=msbm7
\font\fivemsb=msbm5
\textfont\msbfam=\tenmsb
\scriptfont\msbfam=\sevenmsb
\scriptscriptfont\msbfam=\fivemsb

\font\teneufm=eufm10
\font\seveneufm=eufm7
\font\fiveeufm=eufm5
\textfont\eufmfam=\teneufm
\scriptfont\eufmfam=\seveneufm
\scriptscriptfont\eufmfam=\fiveeufm

\def\({\left(}
\def\){\right)}
\def\ee#1{{\vec {\bf e}_{#1}}}
\let\neper=e
\let\ii=i

\def\nep#1{ \neper^{#1}}

\def\tc{\thsp | \thsp}
\let\<=\langle
\let\>=\rangle

\outer\def\nproclaim#1 [#2]#3. #4\par{\medbreak \noindent
   \talato(#2){\bf #1 \Thm[#2]#3.\enspace }%
   {\sl #4\par }\ifdim \lastskip <\medskipamount 
   \removelastskip \penalty 55\medskip \fi}
\def\thmm[#1]{#1}
\def\teo[#1]{#1}
\def\sttilde#1{%
\dimen2=\fontdimen5\textfont0
\setbox0=\hbox{$\mathchar"7E$}
\setbox1=\hbox{$\scriptstyle #1$}
\dimen0=\wd0
\dimen1=\wd1
\advance\dimen1 by -\dimen0
\divide\dimen1 by 2
\vbox{\offinterlineskip%
   \moveright\dimen1 \box0 \kern - \dimen2\box1}
}
\begin{document}

\title[Convex Entropy Decay]
{Convex Entropy Decay via the
  Bochner--Bakry--Emery approach}

\begin{abstract}
 We develop a method, based on a Bochner-type identity, to obtain
 estimates on the exponential rate of decay of the relative entropy
 from equilibrium of Markov processes in discrete settings. 
When this method applies the
 relative entropy decays 
in a convex way. The method is shown to be rather powerful when 
applied to a class of birth and death
processes. 
We then consider other examples, including inhomogeneous zero-range
processes and Bernoulli-Laplace models. For these two models, known
results were limited to the homogeneous case, and obtained via the
martingale approach, whose applicability to inhomogeneous models is
still 
unclear.

%\vskip.45cm

%\noindent
%{\SMALL RESUM\'E}. 
%Nous d\'eveloppons une m\'ethode, insipr\'ee par une
%identit\'e
%de Bochner, pour obtenir des estim\'ees sur la decroissance
%exponentielle
%de l'entropie relative de processus de Markov avec sauts. Lorsque nous
%pouvons aplliquer cette m\'ethode, l'entropie relative est une
%fonction convexe du temps. On montre que la m\'ethode s'applique de
%facon efficace \`a une large classe de processus de naissance et
%mort. On consid\`ere aussi d'autres exemples, comme les processus de
%zero--range et de Bernoulli--Laplace dans des cas
%non--homog\`enes. Pour ces dernier  mod\`eles les r\'esultats connus,
%obtenus par la m\'ethode de martingale, \'etaient limit\'es au cas homog\`ene.

\vskip.5cm



\end{abstract}
\author[P. Caputo]{Pietro Caputo}
\address{Dip. Matematica, Universita' di Roma Tre, L.go S. Murialdo 1,
00146 Roma, Italy} \email{caputo\@@mat.uniroma3.it}
\author[P. Dai Pra]{Paolo Dai Pra}
\address{Dipartimento di Matematica Pura e Applicata,
Universit\`{a} di Padova, Via Belzoni 7, 35131 Padova, Italy}
\email{daipra@math.unipd.it}
\author[G. Posta]{Gustavo Posta} \address{Dip. Matematica, 
Politecnico di Milano, 
P.za Leonardo da Vinci 32, 
I-20133 Milano, Italy} 
\email{gustavo.posta\@@polimi.it}

%\foottext{{\em 2000 MSC: 60K35}}

\date{December 11, 2007}

\maketitle

\thispagestyle{empty}

\section{Introduction}
In the family of functional inequalities that are related to the
convergence to equilibrium of Markov processes (Poincar\`e,
logarithmic Sobolev and Nash inequalities) the {\em modified
  logarithmic Sobolev inequality} ({\bf MLSI}) has been the last to
attract interest among mathematicians, and certainly the less
studied. Consider a time-homogeneous Markov process $(X_t)_{t \geq
  0}$, with values on a measurable space $(S, {\mathcal{S}})$, having
an invariant measure $\pi$. We assume the semigroup $(T_t)_{t \geq 0}$
defined on $L^2(\pi)$ 
by
\[
T_t f(x) := E[f(X_t)|X_0 = x]
\]
is strongly right-continuous, so that the infinitesimal generator
$\cL$ exists, i.e. $T_t = e^{t\cL}$. We also define the nonnegative
quadratic form on $\cD(\cL) \times \cD(\cL)$, called {\em Dirichlet
  form} 
of $\cL$,
\[
\cE(f,g) := - \pi[f \cL g],
\]
where $\cD(\cL)$ is the domain of $\cL$, and we use the notation
$\pi[f]$ for $\int f d\pi$. Given a probability measure $\mu$ on $(S,
{\mathcal{S}})$, we denote by $\mu T_t$ the distribution of $X_t$
assuming $X_0$ is distributed 
according to $\mu$, i.e.
\[
\int f d(\mu T_t) := \int (T_t f) d\mu.
\]
An ergodic Markov process, in particular a countable-state, irreducible and recurrent one, has a unique invariant measure $\pi$, and the rate of convergence of $\mu T_t$ to $\pi$ is a major topic of research. Quantitative estimates on this rate of convergence can be obtained by analyzing functional inequalities. To set up the necessary notations, define the {\em relative entropy} $h(\mu|\pi)$ of the probability $\mu$ with respect to $\pi$ by
\[
h(\mu|\pi) := \pi\left[\frac{d\mu}{d\pi} \log \frac{d\mu}{d\pi} \right],
\]
where  $h(\mu|\pi)$ is meant to be infinite whenever $\mu \not\ll \pi$
or $\frac{d\mu}{d\pi} \log \frac{d\mu}{d\pi} \not\in
L^1(\pi)$. Although $h(\cdot \, | \, \cdot)$ is not a metric in the
usual sense, its use as ``pseudo-distance'' is motivated by a number
of relevant properties, the most 
basic ones being:
\[
h(\mu|\pi) = 0 \ \iff \ \mu = \pi
\]
and
\be{pinsker}
\|\mu - \pi\|_{TV}^2 \leq h(\mu|\pi),
\end{equation}
where $\| \cdot \|_{TV}$ denotes the total variation norm. 
For a generic measurable
function $f\geq 0$ it is common to write
\[
\ent_{\pi}(f) := \left\{ \begin{array}{ll} \pi[f \log f] - \pi[f] \log \pi[f] & \mbox{if } f \log f \in L^1(\pi) \\ +\infty & \mbox{otherwise,} \end{array} \right. 
\]
so that $h(\mu|\pi) = \ent_{\pi} \left( \frac{d\mu}{d\pi} \right)$. Ignoring
technical problems concerning the domains of Dirichlet forms, a simple
formal computation shows that
\begin{equation} \label{p1}
\frac{d}{dt} h(\mu T_t | \pi) = - \cE(T_t^*f, \log T_t^* f)
\end{equation}
where $f := \frac{d\mu}{d\pi}$, $\cL^*$ is the adjoint of $\cL$ in
$L^2(\pi)$, and $T_t^* := e^{t \cL^*}$. 
Therefore,  assuming that, for each $f \geq 0$
\be{MLSI}
\ent_{\pi}(f) \leq \frac{1}{\a}\, \cE(f, \log  f)
\end{equation}
with $\a >0$ (independent of $f$), then (\ref{p1}) can be closed to
get a differential inequality, 
obtaining
\[
 h(\mu T_t | \pi) \leq e^{- \a t}  h(\mu| \pi).
 \]
In other words, estimates on the best constant $\a$ for which the 
functional inequality (\ref{MLSI}) holds provide estimates 
for the rate of exponential convergence to equilibrium of the process, 
in the relative entropy sense.

We shall be interested in {\em reversible} dynamics, i.e.\ 
when $\cL = \cL^*$. When (\ref{MLSI}) holds we say that the 
pair $(\cL,\pi)$ satisfies 
%The estimate (\ref{MLSI}) will
%, it is customary to consider the inequality
%\be{MLSI}
%\ent_{\pi}(f) \leq \frac{1}{\a} \cE(f, \log  f),
%\end{equation}
%which is 
%be called 
the {\em modified logarithmic Sobolev inequality} ({\bf MLSI}) with
constant $\a$. 
This inequality turns out to be intermediate, in a sense that we will 
make precise in a moment, between two more ``traditional'' 
functional inequalities, namely the {\em logarithmic Sobolev inequality} ({\bf LSI})
\be{LSI}
\ent_{\pi}(f^2) \leq \frac{1}{\b}\, \cE(f,   f),
\end{equation}
and the {\em Poincar\'e inequality} ({\bf PI})
\be{PI}
\var_{\pi}(f) \leq  \frac{1}{\g}\, \cE(f,f),
\end{equation}
where $\var_{\pi}(f) := \pi\left[(f-\pi[f])^2 \right]$. It is well
known 
that ({\bf LSI}) is equivalent to hypercontractivity of the 
semigroup $T_t$, i.e. $T_t$ is contractive as linear operator 
from $L^2(\pi)$ and $L^p(\pi)$ for some $p>2$, and ({\bf PI}) 
is equivalent to exponential convergence to equilibrium in $L^2$, 
i.e. $\|T_t f - \pi[f]\|_{L^2(\pi)} \leq e^{-\g t} \| f -
\pi[f]\|_{L^2(\pi)}$. 
Moreover, if we let $\a,\b, \g$ denote the best constant in the 
respective inequality (with the convention that the ``best'' 
constant is zero when the inequality fails for every positive constant), then for reversible systems
\be{comp}
2\g \geq \a \geq 4\b.
\end{equation}
We refer to \cite{DSC} and \cite{BT} for tutorial references on these
inequalities (even though ({\bf MLSI}) is never explicitly mentioned in
\cite{DSC}). It should also be noticed that in the case $\cL$ is the
generator of a reversible diffusion process, e.g.\ $\cL = \frac{1}{2}
\Delta + \nabla V \cdot \nabla$, ({\bf LSI}) and ({\bf MLSI})
coincide. This equivalence, that simply follows from the fact that
$\nabla \log f = \nabla f / f$, does not extend to Markov processes
with jumps; even for a two-state Markov chain, the best constants in
({\bf LSI}) and ({\bf MLSI}) behave quite differently in terms of the
parameter of the invariant measure (see \cite{BobLed}). The case of
processes with jumps leaves some freedom in deciding which inequality 
is the best analogue of the ({\bf LSI}) and there are several
inequalities which are commonly referred to as 
``modified logarithmic Sobolev'' in the literature (all of them
coincide with the ({\bf LSI}) in the diffusion case). Here we only
consider the ({\bf MLSI}) defined in (\ref{MLSI}). Besides the exponential
decay of entropy it is known that this estimate implies useful 
concentration bounds, see \cite{BT}, and thanks to (\ref{pinsker}) it
is a natural tool to estimate mixing times, see \cite{DSC}.
Furthermore,  
we mention that
there is a further family of inequalities interpolating between the
exponential decay in the $L^2$--sense of ({\bf PI}) and the 
exponential decay in the $L\log L$--sense of ({\bf MLSI}) that has
received growing attention in the literature. 
These so--called Beckner
inequalities deal with the exponential decay in the $L^p$--sense, 
for $p\in(1,2)$, see \cite{BT} and references therein for more details.

While the study of ({\bf PI}) and ({\bf LSI}) for large-scale systems
dates 
back to \cite{ma:ol} and \cite{Ya}, 
a similar analysis for ({\bf MLSI}) has been first  
proposed in \cite{Da:Pa:Po}, which deals with Glauber type
dynamics with unbounded particle number per site; for a 
class of such systems ({\bf LSI}) fails, while ({\bf MLSI}) holds with a
strictly positive constant. Dynamics with exchange of particles are 
less understood, with some relevant exceptions (see \cite{GQ,Goel,BT}).

The purpose of this paper is to partially extend to ({\bf MLSI}) the
general tools developed in \cite{bcdp} for obtaining estimates on the
{\em spectral gap}, i.e. the best constant in ({\bf PI}). The results in
\cite{bcdp} apply to interacting particle systems ideas of S. Bochner
(see \cite{boch}), who studied the spectral gap of the Laplacian in
Riemannian manifolds. D.\ Bakry and M.\ Emery \cite{BE} have developed Bochner
ideas, obtaining conditions for ({\bf LSI}) to hold for diffusion
processes. Our work can be interpreted as an attempt to complete the
Bakry and Emery's program for processes with jumps. Although the basic
tools are developed in full generality, useful estimates on the best
constant in ({\bf MLSI}) have been obtained, unfortunately, in a limited
number of examples, if compared both with the diffusion case and the
case of ({\bf PI}) studied in \cite{bcdp}. We believe, however, that
these examples are of interest, and that our approach is 
promising and not yet fully exploited.

\section{The Bochner-Bakry-Emery approach to ({\bf MLSI})} \label{bbe}

We recall here a simple argument relating ({\bf MLSI}) to exponential
decay of entropy and of its time derivative along the semigroup. 
%Below, we always consider the reversible
%case $\cL=\cL^*$. 
To avoid problems that are not relevant
for the applications we have in mind, we assume the state space $S$ of
the Markov chain to be finite or countable. By simple calculus one
checks that, 
for $f>0$
\be{bbe1}
\frac{d}{dt} \ent_{\pi}(T_t f) = - \cE(T_t f, T_t \log f).
\end{equation}
Therefore ({\bf MLSI}) is equivalent to exponential decay of relative
entropy in the sense that for every $\a\geq 0$ one has:
%Thus, if ({\bf MLSI}) holds,
\be{equi}
\a\,\ent_{\pi}(f) \leq \cE(f, \log f)\,, \;\;\; \forall f>0
%\text{for all $f>0$}
\ \ %\Rightarrow 
\iff \ \ 
\ent_{\pi}(T_t f) \leq e^{-\a t} \ent_{\pi}(f)\,,\;\;\;\forall
f>0\,,\;\forall t\geq 0
%\text{for all $f>0$}
\,.
\end{equation}
%\frac{d}{dt} \ent_{\pi}(T_t f) \leq - \a \ent_{pi}(T_t f) \ \ \Rightarrow \ \ \ent_{pi}(T_t f) \leq e^{-\a t} \ent_{pi}(f).
%\end{equation}
Indeed, the implication $\Rightarrow$ is obtained by integrating 
(\ref{bbe1}) and the implication $\Leftarrow$ follows by subtracting
$\ent_{\pi}(f)$ from the right hand side of (\ref{equi}), dividing by
$t$, and taking $t\to 0 $.

In other words, ({\bf MLSI}) is equivalent to a control on the first
time derivative of entropy. The following simple Lemma 
(which is well known, see e.g. \cite{Le}) is based on a similar 
control of second derivatives.
\bl{l1}
Suppose the generator  $\cL$ is self-adjoint in $L^2(\pi)$ and the
resulting Markov chain is irreducible. % , and that the inequality
Then, for every $\k \geq 0$ we have the equivalence:
\begin{align}
\la{MLSI''}
&\k \,
\cE(f, \log f) \leq \pi[\cL f \cL\log f] + \pi\left[ \frac{(\cL f)^2}{f} \right]\,,\;\;\;\forall f>0%\text{for all $f>0$}
 \\ 
&\quad\quad \quad \iff \ \ 
\cE(T_t f,\log T_t f) \leq e^{-\k t} \cE(f,\log f)\,,
\;\;\;\forall f>0\,,\;\forall t\geq 0%\text{for all $f>0$}
\,.\nonumber
\end{align}
%\be{MLSI'}
%k \cE(f \log f) \leq \pi[\cL f \cL\log f] + \pi\left[ \frac{(\cL f)^2}{f} \right]
%\end{equation}
%holds for every $f>0$. 
Moreover, if (\ref{MLSI''}) holds for some $\k$, then 
{\rm ({\bf MLSI})} holds with $\a = \k$.
%Then ({\bf MLSI}) holds with $\a = \k$.
\el
\bpr
Computing second derivatives we obtain, 
%also using (\ref{bbe1}) and (\ref{MLSI'})
\begin{equation}
\frac{d^{2}}{dt^{2}} \ent_{\pi}(T_t f) = - 
\frac{d}{dt} \cE(T_t f,\log T_t f) = 
\pi\left[ \cL T_t f \cL \log T_t f\right] + 
\pi \left[ \frac{(\cL T_t f)^2}{T_t f} \right]\,.
\end{equation}
%\begin{align*}
%\frac{d^2}{dt^2} Ent_{\pi}(S_t f) &= - \frac{d}{dt} \cE(T_t f,\log T_t f) = \pi\left[ \cL T_t f \cL \log T_t f\right] + \pi \left[ \frac{(\cL T_t f)^2}{T_t f} \right]  \\ & \geq \k \cE(T_t f,\log T_t f) = -\k\frac{d}{dt} \ent_{\pi}(T_t f).
%\end{align*}
The equivalence (\ref{MLSI''}) therefore follows as in the case of
(\ref{equi}) discussed above. 
%which implies
%\[
%0 \leq \cE(T_t f,\log T_t f) \leq e^{-\k t} \cE(f,\log f),
%\]
%where the nonnegativity comes from the well known inequality  $ \cE(f,\log f) \\ \geq  4 \cE(\sqrt{f},\sqrt{f}) \geq 0$ (see \cite{DSC}). Thus the inequality
To prove the last assertion, note that the inequality
\[
\k \,\ent_{\pi}(f) \leq  \cE(f,\log f)\,,
\]
follows by integrating from $0$ to $\infty$ the inequality
\[
- \frac{d}{dt} \cE(T_t f,\log T_t f) \geq -\k\frac{d}{dt} \ent_{\pi}(T_t f)\,.
\]
\epr
\br{r1}
Inequality (\ref{MLSI''}) is in general {\em strictly} stronger than
({\bf MLSI}): it implies uniform exponential decay of entropy, but also that
the decay is convex in time. While it is easily seen that 
%well known inequalities (see e.g.\
%\cite{DSC}) imply 
$\cE(f,\log f) \geq  0$ for all $f>0$, %4 \cE(\sqrt{f},\sqrt{f})
%\geq 0$, 
nothing forces the second derivative of $\ent_{\pi}(T_t f)$
to be non--negative.  
There are examples showing that ({\bf MLSI}) may hold without convexity
in time of entropy. An example in the continuous setting, due to
B.\ Helffer, can be found in \cite{Le}. 
An example in the discrete setting will be given later in this paper,
see Section 4.
\er

In order to investigate the validity of (\ref{MLSI''}) in the discrete
setting, we write the generators of our Markov chains in the form
\be{geno}
\cL f(\eta ) = \sum_{\g \in G} c(\eta ,\g) [f(\g(\eta )) - f(\eta ) ] =: \sum_{\g \in G} c(\eta ,\g) \nabla_{\g} f (\eta )
\end{equation}
where $G$ is some (finite or countable) set of functions from $S$ to
$S$ (the {\em allowed moves}) and $c:S\times G\to [0,\infty)$
represent the jump rates. It is easily seen that the generator of
every finite or countable Markov chains can be written in this form; a
form that, as we will see, becomes rather natural in many specific
examples. 

With these notations, reversibility is expressed as 
follows. 
\begin{itemize}
\item[({\bf Rev})] For every $\g \in G$ there exists $\g^{-1} \in G$ such that $\g^{-1} \g(\eta ) = \eta $ for every $\eta  \in S$ such that $c(\eta ,\g)>0$. Moreover for every $f:S \ra \R$ bounded
%$F:S \times G \ra \R$ bounded
\[
\pi \left[ c(\eta ,\g) f(\eta ) \right] = \pi \left[ c(\eta ,\g^{-1})
  f(\g^{-1}(\eta )) \right]. %\what
\]
\end{itemize}
Under ({\bf Rev}) it is easy to see that
\[
\cE(f,g)  = \frac{1}{2}\, \pi \left[ \sum_{\g \in G} c(\eta ,\g) \nabla_{\g} f (\eta ) \nabla_{\g} g(\eta ) \right].
\]
In particular the Dirichlet form is symmetric, so $\cL$ is self-adjoint in $L^2(\pi)$.

One of the key point of the Bochner-Bakry-Emery approach is the so called {\em Bochner identity}; a version of this identity in the discrete setting is given in next Lemma.
\bl{l2}
Let $R: S \times G \times G \ra [0,+\infty)$ be such that 
\begin{eqnarray*}
({\bf P1}): &  R(\eta ,\g,\d) = R(\eta ,\d,\g) & \mbox{$\forall \;\eta ,\g,\d$ with $R(\eta ,\g,\d)>0$} \\ 
({\bf P2}): &   \pi \left[ \sum_{\g,\d}
R(\eta ,\g,\d) \psi(\eta ,\g,\d) \right]  = \pi \left[  
\sum_{\g,\d}R(\eta ,\g,\d) \psi(\g(\eta ), \g^{-1},\d) \right]&
\mbox{$\forall \;\psi:S \times G \times G$ bounded} \;%\what
\\
({\bf P3}): &  \g \d (\eta ) = \d \g (\eta ) & \mbox{$\forall \;\eta ,\g,\d$ with $R(\eta ,\g,\d)>0$}
\end{eqnarray*}
Then, for every $f,g$ the following {\em  Bochner}-type  identity holds
\[
\pi\left[ \sum_{\g,\d} R(\eta ,\g,\d) \nabla_{\g} f(\eta )
  \nabla_{\d}\,g(\eta ) \right]  = \frac{1}{4}\, \pi \left[ \sum_{\g,\d}
  R(\eta ,\g,\d) 
\nabla_{\g}\nabla_{\d} f(\eta ) \nabla_{\g}\nabla_{\d}\,g(\eta ) \right] 
\]
\el
\bpr
First, by ({\bf P3}), $\nabla_{\g}\nabla_{\d} f(\eta ) \nabla_{\g}\nabla_{\d}g(\eta ) = \nabla_{\g}\nabla_{\d} f(\eta ) \nabla_{\d}\nabla_{\g}g(\eta )$ whenever $R(\eta ,\g,\d)>0$. Then
write
\[
 \nabla_{\g}\nabla_{\d} f(\eta ) \nabla_{\d}\nabla_{\g}g(\eta ) = \nabla_{\d}f(\g(\eta ))\nabla_{\g}g(\d(\eta )) - \nabla_{\d} f(\g(\eta )) \nabla_{\g}g(\eta )  - \nabla_{\d} f(\eta ) \nabla_{\g}g(\d(\eta )) + \nabla_{\d} f(\eta ) \nabla_{\g}g(\eta )
 \]
 We show that each one of the four summands in the r.h.s. of this last formula, when multiplied by $R(\eta ,\g,\d)$, summed over $\g,\d$ and averaged over $\pi$ gives
 \[
\pi\left[ \sum_{\g,\d} R(\eta ,\g,\d) \nabla_{\d} f(\eta ) \nabla_{\g}g(\eta )\right].
 \]
 For the fourth summand there is nothing to prove. Moreover, by ({\bf P2}),
 \begin{align*}
\pi\left[ \sum_{\g,\d} R(\eta ,\g,\d) \nabla_{\d} f(\eta )
  \nabla_{\g}g(\eta )\right]  &= \pi\left[  \sum_{\g,\d}R(\eta ,\g,\d)
  \nabla_{\d}f(\g(\eta )) \nabla_{\g^{-1}} g(\g(\eta )) \right] \\ &= -
\pi\left[  \sum_{\g,\d}R(\eta ,\g,\d) 
\nabla_{\d}f(\g(\eta )) \nabla_{\g} g(\eta ) \right]
\end{align*}
which takes care of the second and, by symmetry, of the third summand. For the first summand we use first ({\bf P2}), then ({\bf P1}), ({\bf P2}) again and ({\bf P3}):
\begin{align*}
&\pi\left[ \sum_{\g,\d} R(\eta ,\g,\d) \nabla_{\d} f(\eta )
  \nabla_{\g}g(\eta )\right]  = \pi\left[  \sum_{\g,\d}R(\eta ,\g,\d)
  \nabla_{\d}f(\g(\eta )) \nabla_{\g^{-1}} g(\g(\eta )) \right] \\ &\qquad = -
\pi\left[ \sum_{\g,\d} R(\eta ,\g,\d) \nabla_{\g} f(\d(\eta ))
  \nabla_{\d}g(\eta ) \right] = - \pi\left[ \sum_{\g,\d} R(\eta
  ,\g,\d) \nabla_{\g^{-1}} f(\d\g(\eta )) \nabla_{\d}g(\g(\eta ))
\right] \\ &\qquad = \pi\left[ \sum_{\g,\d} R(\eta ,\g,\d) \nabla_{\g}
  f(\d(\eta )) 
\nabla_{\d}g(\g(\eta )) \right] 
\end{align*}
\epr

\bc{c1}
Let $R: S \times G \times G \ra [0,+\infty)$ be such that {\rm ({\bf P1})}, {\rm ({\bf P2})} and {\rm ({\bf P3})} hold.
Define
\[
\G(\eta ,\g,\d) := c(\eta ,\g) c(\eta ,\d) - R(\eta ,\g,\d).
\]
Then, for every $f>0$:
 \[
\pi[\cL f \cL\log f] + \pi\left[ \frac{(\cL f)^2}{f} \right] \\ \geq \pi\left[ \sum_{\g,\d}\Gamma(\eta ,\g,\d)\left( \nabla_{\g} f(\eta ) \nabla_{\d} \log f(\eta ) +  \frac{\nabla_{\g} f (\eta ) \nabla_{\d} f(\eta )}{f(\eta )}\right) \right].
\]
\ec
\bpr
First observe that
\begin{align}\label{bbe2}
\pi[\cL f \cL\log f] + \pi\left[ \frac{(\cL f)^2}{f} \right] 
& = \pi\left[ \sum_{\g,\d}c(\eta ,\g)c(\eta ,\d) \nabla_{\g} f(\eta )
  \nabla_{\d} \log f(\eta ) \right] \nonumber \\ &\quad
 + \pi\left[ \sum_{\g,\d}c(\eta
  ,\g)c(\eta ,\d)\frac{\nabla_{\g}f(\eta ) \nabla_{\d} f(\eta
    )}{f(\eta )} 
\right].
\end{align}
Now, if we apply Bochner's identity to the first summand of the right hand side of (\ref{bbe2}) we obtain
\begin{multline}\label{bbe3}
\pi[\cL f \cL\log f] + \pi\left[ \frac{(\cL f)^2}{f} \right] 
 = \pi\left[ \sum_{\g,\d}\Gamma (\eta ,\g,\d) \nabla_{\g} f(\eta ) \nabla_{\d} \log f(\eta ) \right] \\ + \frac{1}{4}\,\pi\left[ \sum_{\g,\d}R(\eta ,\g,\d) \nabla_{\g} \nabla_{\d}f(\eta ) \nabla_{\g}\nabla_{\d} \log f(\eta ) \right] + \pi\left[ \sum_{\g,\d}c(\eta ,\g)c(\eta ,\d)\frac{\nabla_{\g}f(\eta ) \nabla_{\d} f(\eta )}{f(\eta )} \right].
\end{multline}
Thus the conclusion follows immediately if we show that
\begin{align}\label{bbe4}
\frac{1}{4}\, &\pi\left[ \sum_{\g,\d}R(\eta ,\g,\d) \nabla_{\g} \nabla_{\d}f(\eta ) \nabla_{\g}\nabla_{\d} \log f(\eta ) \right]  + \pi\left[ \sum_{\g,\d}c(\eta ,\g)c(\eta ,\d)\frac{\nabla_{\g}f(\eta ) \nabla_{\d} f(\eta )}{f(\eta )} \right] \nonumber\\ &\quad\quad\geq 
 \pi\left[ \sum_{\g,\d}\Gamma(\eta ,\g,\d) \frac{\nabla_{\g} f (\eta ) \nabla_{\d} f(\eta )}{f(\eta )} \right]\,,
 \end{align}
 or, equivalently
 \be{bbe5}
 \frac{1}{4}\, \pi\left[ \sum_{\g,\d}R(\eta ,\g,\d) \nabla_{\g} \nabla_{\d}f(\eta ) \nabla_{\g}\nabla_{\d} \log f(\eta ) \right] \\ + \pi\left[ \sum_{\g,\d}R(\eta ,\g,\d)\frac{\nabla_{\g}f(\eta ) \nabla_{\d} f(\eta )}{f(\eta )} \right] \geq 0.
 \end{equation}
 Now we apply to the second summand in (\ref{bbe5}) the same argument used in the proof of Lemma \ref{l2}. The result is more cumbersome, due to the presence of the denominator $f(\eta )$:
 \begin{align*}
 &\pi\left[ \sum_{\g,\d}R(\eta ,\g,\d)\frac{\nabla_{\g}f(\eta ) \nabla_{\d} f(\eta )}{f(\eta )} \right] \\ &\qquad\quad
 =  \frac{1}{4}\,\pi\left[ \sum_{\g,\d}R(\eta ,\g,\d) \left\{\nabla_{\g}\left(\frac{\nabla_{\d} f(\eta )}{f(\d(\eta ))} \right) \nabla_{\g} \nabla_{\d} f(\eta ) - \nabla_{\g} \left( \frac{(\nabla_{\d} f(\eta ))^2}{f(\eta )f(\d(\eta ))} \right) \nabla_{\g} f(\eta ) \right\} \right].
 \end{align*}
 Substituting in (\ref{bbe5}) we get
 \begin{align} \label{bbe6}
  &\frac{1}{4}\, \pi\left[ \sum_{\g,\d}R(\eta ,\g,\d) \nabla_{\g}
    \nabla_{\d}f(\eta ) \nabla_{\g}\nabla_{\d} \log f(\eta ) \right] +
  \pi\left[ \sum_{\g,\d}R(\eta ,\g,\d)\frac{\nabla_{\g}f(\eta )
      \nabla_{\d} f(\eta )}{f(\eta )} \right]\nonumber \\& \qquad\quad =
  \pi\Bigg[
    \sum_{\g,\d}R(\eta ,\g,\d) \bigg\{ 
\nabla_{\g} \nabla_{\d}f(\eta )
      \nabla_{\g}\nabla_{\d} \log f(\eta ) 
\nonumber \\ & \qquad\qquad\qquad
  + \nabla_{\g}\left(\frac{\nabla_{\d} f(\eta )}{f(\d(\eta
        ))} \right) \nabla_{\g} \nabla_{\d} f(\eta ) - \nabla_{\g}
    \left( \frac{(\nabla_{\d} f(\eta ))^2}{f(\eta )f(\d(\eta ))}
    \right) \nabla_{\g} f(\eta ) \bigg\} 
\Bigg] \,.
\end{align}
Setting $a:=f(\eta ), b:=f(\d(\eta )), c := f(\g(\eta )), d := f(\d \g(\eta ))$, one checks that the term in braces in the right hand side of (\ref{bbe6}) equals the sum of the following $4$ expressions
\[
\begin{array}{c}
d\log d - d \log(bc/a) + (bc/a) - d \\ c\log c - c \log(da/b) + (da/b) -c \\ b\log b - b \log(da/c) + (da/c) - b \\ a \log a - a \log(bc/d) + (bc/d) -a
\end{array}
\]
which are all nonnegative, since $\a \log \a - \a \log \b + \b - \a
\geq 0$ for every $\a,\b >0$. This shows that (\ref{bbe6}) is
nonnegative, which completes the proof. \epr

Summing up, we are led to the following result.
\bp{pr1}
Suppose there exists a constant $\k >0$ such that for every $f>0$
\[
\pi\left[ \sum_{\g,\d}\Gamma(\eta ,\g,\d)\left( \nabla_{\g} f(\eta )
    \nabla_{\d} \log f(\eta ) +  \frac{\nabla_{\g} f (\eta )
      \nabla_{\d} f(\eta )}{f(\eta )}\right) \right]  \geq
\frac{\k}{2}\,  \pi \left[ \sum_{\g \in G} c(\eta ,\g) \nabla_{\g} f
  (\eta ) \nabla_{\g} \log f(\eta ) 
\right].
\]
Then (\ref{MLSI''}) holds. In particular, {\rm ({\bf MLSI})} holds with $\a = \k$.
\ep

We do not have a general choice for the function $R$ of Lemma \ref{l2} to use for applying the criterion in Proposition \ref{p1}. One option is given in \cite{bcdp}, Proposition 2.4. That choice, however, does not work for the examples in Section 4 below. We only mention that in all examples we obtain $R(\eta ,\g,\d)$ by ``modifying'' $c(\eta ,\g) c(\g(\eta ),\d)$ in order to fulfill properties ({\bf P1})--({\bf P3}).

\section{A warming-up example: birth and death processes}

Consider a birth and death process on $\N$ with generator
\be{genbd}
\cL f(n) = a(n) \nabla_+ f(n) + b(n) \nabla_- f(n)\,.
\end{equation}
In the language introduced in the previous section, $G = \{+,-\}$,
where $+(n) = n+1$ and $-(n) = (n-1) {\bf 1}_{n>0}$, and one is the
inverse of the other. In particular, $\nabla_\pm f(n) = f(n\pm 1)-f(n)$. 
The rates $a,b$ are non--negative functions on $\bbN$
such that $b(0)=0$ 
and we assume that there exists a probability $\pi$ on $\N$ 
such that the detailed balance equation 
\be{rever}
a(n) \pi(n) = b(n+1) \pi(n+1)
\end{equation}
holds true. Moreover, we assume that the resulting Markov chain is
irreducible. 

Setting $c(n,+) = a(n)$, $c(n,-) = b(n)$ we see that condition ({\bf
Rev}) in the previous section is satisfied and 
\[
{\cE}(f,g) = \pi\left[ a(n) \nabla_+ f(n) \nabla_+ g(n) \right] = \pi \left[ b(n) \nabla_- f(n) \nabla_- g(n) \right]\,.
\] 
We define $R$ as follows:
\begin{eqnarray}
R(n,+,+) & := & a(n)a(n+1) \nonumber
\\ R(n,-,-) & := & b(n)b(n-1) \label{rdb}
\\ R(n,+,-) = R(n,-,+) & := & a(n) b(n)\,.\nonumber
\end{eqnarray}
It is a simple exercise to show that conditions ({\bf P1})--({\bf P3})
of Lemma \ref{l2}  are satisfied. In particular, ({\bf P2}) follows
by application of reversibility.  Then, letting as before
$\Gamma(n,\d,\g) = c(n,\g)c(n,\d) - R(n,\g,\d)$, Corollary \ref{c1}
yields
%we have 
 \begin{align}
 &\pi[\cL f \cL\log f] + \pi \left[ \frac{(\cL f)^2}{f} \right]   \geq
 \pi \left[ \sum_{\g,\d \in G}\Gamma(n,\d,\g)\left(\nabla_{\g}f(n)
     \nabla_{\d} \log f(n) + \frac{\nabla_{\g} f(n) \nabla_{\d}
       f(n)}{f(n)}\right) \right]\nonumber\\ &\qquad = 
\pi\Big[a(n)[a(n) - a(n+1)] \nabla_+ f(n) \nabla_+ \log f(n) 
 + b(n)[b(n) - b(n-1)] \nabla_- f(n) \nabla_- \log f(n) \Big]
 \nonumber \\ & \qquad\qquad + 
  \pi\left[a(n)[a(n) - a(n+1)] \frac{\left(\nabla_+ f(n)\right)^2}{f(n)} 
%  \right. \\  \left.
  + b(n)[b(n) - b(n-1)] \frac{\left(\nabla_- f(n)\right)^2}{f(n)}
\right] \,.
\la{ibn}
 \end{align}
 
%   
% \geq \pi\left[a(n)[a(n) - a(n+1)] \nabla_+ f(n) \nabla_+ \log f(n)
%   + b(n)[b(n) - b(n-1)] \nabla_- f(n) \nabla_- \log f(n) \right].
We consider the following assumption:
\begin{itemize}
\item[({\bf A})]
$a(n) \geq a(n+1)$ and $b(n+1) \geq b(n)$, and 
there exists $c>0$ 
such that for every $n\geq 0$, \be{aa}
a(n) - a(n+1) + b(n+1) - b(n) \geq c\,.
\end{equation}
\end{itemize}

Assuming ({\bf A}), (\ref{ibn}) can be further estimated as follows. 
Thanks to monotonicity of the rates we can drop the terms in the last
line of (\ref{ibn}). Moreover from the reversibility (\ref{rever}) we
see that 
$$
\pi\Big[b(n)[b(n) - b(n-1)] \nabla_- f(n) \nabla_- \log f(n) \Big]
= \pi\Big[a(n)[b(n+1) - b(n)] \nabla_+ f(n) \nabla_+ \log f(n) \Big]\,.
$$
Therefore, using (\ref{aa}) we arrive at 
$$
\pi[\cL f \cL\log f] + \pi \left[ \frac{(\cL f)^2}{f} \right]   \geq
c\,\pi\left[a(n) \nabla_+ f(n) \nabla_+ \log f(n) \right] =
c\,\cE(f,\log f)\,.
$$
Recalling Proposition \ref{pr1} we have therefore proved the following result.
\begin{theorem}\label{t1}
Under assumption {\rm ({\bf A})}, both (\ref{MLSI''}) and {\rm ({\bf MLSI})} hold 
with constant $c$.
\end{theorem}
%\subsection{Some applications of Theorem \ref{t1}}
There are well known criteria for the validity of ({\bf PI})
or ({\bf LSI}) for one-dimensional processes as the ones considered
above, see \cite{miclo} for explicit estimates on the constants involved. 
On the other hand, we are not aware of any such result concerning ({\bf MLSI}).
As far as we know Theorem \ref{t1} is the first general 
sufficient condition for the validity of ({\bf MLSI}). 
Moreover, despite of its simplicity, this result is sharper than it
may appear, as the following examples illustrate.
 
\subsection{Poisson case}
The Poisson case refers to the choice $a(n)=\l$, $b(n)=n$ and
$\pi_\l(n)=\frac{\l^n}{n!}\,\nep{-\l}$, with $\l>0$. It is well known
that ({\bf LSI}) fails in this case. To see this, take $f_k(n) := 
{\bf 1}_{(k,+\infty)}(n)$ in (\ref{LSI}) and then let $k\to\infty$. On
the other hand assumption ({\bf A})
is satisfied with $c=1$ so that Theorem \ref{t1} yields the following
estimate
for any $f>0$:
\be{poi}
\ent_{\pi_\l}(f)\leq \l\,\pi_\l\left[\nabla_+ f    \nabla_+ \log f \right]\,.
%\cE(f,\log f)\,.
\end{equation}
Note that this estimate is sharp, in the sense that no better constant
than $\l$ can satisfy (\ref{poi}) for all $f>0$. To see this, it
suffices to take $f_k(n) = \nep{kn}$, for fixed $k\in\bbN$; 
simple computations show that $\ent_{\pi_\l}(f) =
(k\l\nep{k}-\l\nep{k}+\l) 
\nep{\l(\nep{k}-1)}$
while $\pi_\l\left[\nabla_+ f    \nabla_+ \log f \right] =
k(\nep{k}-1) \nep{\l(\nep{k}-1)}$, and equality is approached in
(\ref{poi}) as 
$k\to\infty$.
% 
% up to a factor 2. This can be checked
%by taking $f_k(n) = \nep{-n/k}$, in which case, neglecting $O(k^{-3})$
%terms one has 
%$$\ent_{\pi_\l}(f) = \frac{\l\exp({\l(\nep{-k^{-1}}-1)})}{2k^2}\,,\quad
%\cE(f,\log f) = \frac{\l\exp({\l(\nep{-k^{-1}}-1)})}{k^2}\,.$$
%That is, if $\a$ is the best constant in the ({\bf MLSI}), we must have
%$1\leq \a\leq 2$. 
In the Poisson case one can obtain inequality (\ref{poi}) also 
using the Poisson limit of the binomial distribution as
in \cite{BobLed,Da:Pa:Po}. In this special case
the analysis can be pushed beyond these statements, see e.g.\ \cite{Wu,Chafai}
for further developments.

\subsection{Log--concave probabilities}
A non--negative function 
$\g$ on $\bbN$ is called log--concave if %$\g(n):=\mu(n)$ satisfies
\begin{equation}
\g (n)^2\geq \g(n+1)\,\g(n-1)\,.
\la{LC1}
\end{equation}
%The Poisson measure $\pi_\l$ satisfies (\ref{LC1}) for any $\l>0$.
Suppose our measure $\pi$ is such that $\g(n):=n!\pi(n)$  
satisfies (\ref{LC1}). 
Such a measure is sometimes called ultra log--concave.
If we set $a(n)=1$ for all $n\geq 0$, then it follows that 
$$
b(n+1)-b(n) = \frac{\g(n)}{\g(n+1)}\,(n+1) - \frac{\g(n-1)}{\g(n)}\,n
\geq \frac{\g(n-1)}{\g(n)} \geq \cdots \geq \frac{\g(0)}{\g(1)} = b(1)\,.
$$ 
From Theorem \ref{t1} we obtain that ({\bf MLSI}) holds in this model
with $\a=b(1)$, i.e.\ that $\pi$ satisfies
\be{ulc}
b(1)\,\ent_{\pi}(f) \leq %\frac1{b(1)}\,
\pi\left[\nabla_+ f \nabla_+ \log f \right]\,,
\end{equation}
for any $f>0$. Note that (\ref{poi}) is a special case of (\ref{ulc}).
It has been shown in \cite{Johnson} that Poisson measures maximize
entropy in the class of ultra log--concave measures. It is interesting
to note that the convexity results obtained in \cite{Johnson} can be
derived in a simple way from the arguments in our proof of
Theorem \ref{t1}.

\subsection{Random walks}\label{rws}
Another example is the simple random walk on a segment
$[0,N]\cap\bbZ$ with reflecting boundary conditions. 
Here the ({\bf MLSI}) constant is
known to be of order $1/N^2$ by (\ref{comp}), 
since both ({\bf LSI}) and ({\bf PI}) 
can be shown to hold with $\b\sim 1/N^2$ and $\g\sim 1/N^2$
(see e.g.\ \cite{Ya} for the proof of the ({\bf LSI})). 
Let us show that this can be deduced from the
above bounds. Let $\mu$ denote the uniform probability over
$[0,N]\cap\bbZ$. We want to prove that for some constant $C$, for all
$f>0$:
\be{flatrw}
\ent_{\mu}(f) \leq C\,N^2\,%\frac1{b(1)}\,
\mu\left[\nabla_+ f \nabla_+ \log f \right]\,.
\end{equation}
Let $\pi$ denote the probability on $[0,N]\cap\bbZ$ such that
$\pi(n)$ is proportional to ${\bf 1}_{\{0\leq n \leq N\}}\,\nep{-n^2/N^2}$. 
It is easy to check that $\pi$ is
equivalent to $\mu$, i.e.\ 
$\d\leq \mu(n)/\pi(n)\leq \d^{-1}$ for some $\d>0$ 
independent of $N$ and $n$. Then, by a standard comparison argument
(see e.g.\ Lemma 3.3 in \cite{DSC}), it is
sufficient to prove (\ref{flatrw}) for the measure $\pi$ in place of
$\mu$. This in turn follows from Theorem \ref{t1}. Indeed, setting
$a(n)=1$, for all $0\leq n\leq N-1$ and $a(N)=0$, 
we have $b(n) = \pi(n-1)/\pi(n)$ for
all $1\leq n\leq N$, and $b(0)=0$. Therefore (\ref{aa}) applies with
$c^{-1}=O(N^2)$, for all $0\leq n\leq N-1$ and Theorem \ref{t1} allows to prove the claim.
%
%\be{entdisson}
%\ent_\pi(f)\leq \d^{-1}\,N^2\,\pi\left[\grad_+ f \,\grad_+\log f\right]\,,\quad \; f>0\,.
%\end{equation}

\subsection{Non--monotone rates}
 By means of a perturbation argument we can relax the monotonicity
 requirement in assumption ({\bf A}). More precisely, suppose that
 $a(n)=1$ for all $n\geq 0$, so that the probability measure $\pi$ satisfies
 $\pi(n+1)/\pi(n)=1/b(n+1)$, or
\be{zr}
\pi(n)=\frac{\pi(0)}{b(1)\cdots b(n)}\,,\quad \;n\geq 1\,.
\end{equation}
In this case Theorem \ref{t1} shows that if $b(n+1)-b(n)\geq c$ for
all $n\geq 0$ then ({\bf MLSI}) holds with $\a=c$. The next result,
which is a key ingredient in the proof of the main theorem
in \cite{CP}, shows
that if we impose a Lipschitz condition on the rates $b(n)$, then it is
sufficient to have monotonicity on a large scale.  
\bp{nonmon}
Suppose that there exist
$C_1<\infty$, $\d>0$ and $n_0\in\bbN$ such that 
\begin{gather*}
\sup_{n\geq 0} |b(n+1)-b(n)|\leq C_1\,,\quad \;\text{ and}\quad\;
\inf_{n\geq 0} \left[b(n+n_0)-b(n)\right]\geq \d\,.
\end{gather*}
Then, for some constant $C$ which may depend on $C_1$, $\d$, and $n_0$ only, the probability measure (\ref{zr}) satisfies 
\be{entdisson}
\ent_\pi(f)\leq C\,\pi\left[\grad_+ f \,\grad_+\log f\right]\,,\quad \; f>0\,.
\end{equation}
\ep
To prove the proposition we shall need a preliminary lemma.
Define
\be{ctilde}
\tilde b(k) := b(k) + \frac1{n_0}\sum_{j=1}^{n_0-1}
\frac{n_0-j}{n_0}\,[b(k+j) + b(k-j) - 2b(k)]\,,\quad \;k\geq n_0\,,
\end{equation}
and, when $0\leq k<n_0$, set $\tilde b(k)= \tilde b(n_0)k/{n_0}$.
Let us call $\tilde \pi$ the 
probability measure obtained from $\tilde b$ by (\ref{zr}).
\bl{equiv}
The rate function $\tilde b$ is uniformly increasing:\ there exists $\d_1>0$ such that  
$\grad_+\tilde b \geq \d_1$.
Moreover, $\pi$ and $\tilde\pi$ are equivalent: 
there exists $C>0$ such that $C^{-1}
\leq \tilde \pi(n)/\pi(n)%\frac{\tilde \pi(n)}{\pi(n)}
\leq C$, 
for all $n\in\bbN$.
%\be{equi1}
%\frac1C \,\leq\, \frac{\tilde \mu(n)}{\mu(n)}\, \leq\, C\,.
%\end{equation}
\el
\proof
%To establish (\ref{equi0}) w
We rewrite $\tilde b(k)$, $k\geq n_0$:
\begin{align*}
\tilde b(k) &= \frac{b(k)}{n_0} + \frac1{n_0}\sum_{j=1}^{n_0-1}
\frac{n_0-j}{n_0}\,[b(k+j) + 
b(k-j)] \nonumber\\
& = \frac1{n_0}\sum_{j=0}^{n_0-1}
\left\{\frac{n_0-j}{n_0}\,b(k+j) +  \frac{j}{n_0}\,b(k+j-n_0)\right\}
\,.
\end{align*}
To compute $\grad_+\tilde b$ we use summation by parts in the form
\be{sumpart}
\sum_{j=\ell}^m \psi(j)\grad_+\varphi(j) = 
\psi(m)\varphi(m+1) - \psi(\ell)\varphi(\ell)
- \sum_{j=\ell+1}^m \varphi(j)\grad_+\psi(j-1)\,,
\end{equation}
where $\ell<m$ and $\psi,\varphi$ are arbitrary functions. We apply 
(\ref{sumpart}) with $\ell=0,m=n_0-1$, first to the case 
$\psi(j)=(n_0-j)/n_0$, $\varphi(j)=b(k+j)$ and then to the case
$\psi(j)={j}/{n_0}$, $\varphi(j)=b(k+j-n_0)$. The conclusion is that,
for every $k\geq n_0$ we have
\be{ctilde3}
\grad_+\tilde b(k) = \frac1{n_0^2}\sum_{j=0}^{n_0-1}
\,[b(k+j)  -
b(k+j-n_0)]\,.
\end{equation}
Since 
$\grad_+ \tilde b(k)\geq \tilde b(n_0)/n_0$ for every $k<n_0$, the 
claim $\grad_+ \tilde b\geq \d_1$ follows 
from (\ref{ctilde3}) and 
the hypothesis $b(n+n_0)  - b(n) \geq \d$.

We turn to the proof of the equivalence of $\pi,\tilde\pi$. 
We have to prove that there exists 
$C\in[1,\infty)$ such that for every $n\in\bbN$
$$%{equi4}
C^{-1} \,\leq\, \prod_{k=1}^n\frac{\tilde b(k)}{b(k)}\, \leq\, C\,.
$$%\end{equation}
We shall 
prove the left inequality above. The right inequality is obtained with 
the same proof by interchanging the role of $b$ and $\tilde b$.
Passing to logarithms it suffices to prove
\be{series}
\sup_n\sum_{k=1}^n \frac{b(k)-\tilde b(k)}{\tilde b(k)} < \infty\,.
\end{equation}
From (\ref{ctilde}), writing
$$
b(k+j) + b(k-j) - 2b(k)
= \sum_{i=0}^{j-1} \left[\grad_+ b(k+i)-\grad_+ b(k+i-j)\right]\,,
$$
we have
\be{equi5}
\sum_{k=n_0}^n \frac{b(k)-\tilde b(k)}{\tilde b(k)}
= \frac1{n_0}\sum_{j=1}^{n_0-1}\frac{n_0-j}{n_0}\sum_{i=0}^{j-1}
\sum_{k=n_0}^n 
\frac{\grad_+ b(k+i-j)-\grad_+ b(k+i)}{\tilde b(k)}
\end{equation}
Now, for every fixed $i<j$ 
we can use summation by parts as in (\ref{sumpart}), with $\ell=n_0,m=n$
and $\psi(k)=1/\tilde b(k)$, $\varphi(k)=b(k+i-j)$ to obtain
\begin{align*}
\sum_{k=n_0}^n & 
\frac{\grad_+ b(k+i-j)}{\tilde b(k)} \\ & \;=\,
\frac{b(n+1+i-j)}{\tilde b(n)} - \frac{b(n_0+i-j)}{\tilde b(n_0)}
+ \sum_{k=n_0+1}^n \frac{b(k+i-j)\grad_+\tilde b(k-1)}{\tilde b(k)\tilde b(k-1)}\,.
\end{align*}
Another application of (\ref{sumpart}) with $\varphi(k)=b(k+i)$ yields
therefore the identity
\begin{align}
\sum_{k=n_0}^n 
&\frac{\grad_+ b(k+i-j)-\grad_+ b(k+i)}{\tilde b(k)} 
 \label{equi6}\\
& = 
\frac{b(n+1+i-j)-b(n+1+i)}{\tilde b(n)} - 
\frac{b(n_0+i-j)-b(n_0+i)}{\tilde b(n_0)}\, +\, \nonumber\\
& \quad + \sum_{k=n_0+1}^n \frac{[b(k+i-j)-b(k+i)]
\grad_+\tilde b(k-1)}{\tilde b(k)\tilde b(k-1)}
\nonumber
\end{align}
Since $\tilde b(k)\geq \d_1 k$, 
the sequence $(\tilde b(k)\tilde b(k-1))^{-1}$ is 
summable. By hypothesis the increments of $b$ (and $\tilde b$) 
are uniformly bounded 
and therefore the sum in (\ref{equi6}) is uniformly bounded in $n$, 
for every $i<j<n_0$. Now (\ref{series}) follows from (\ref{equi5}).
\qed
%
%
% 
%\begin{proposition}\label{pr2}
%Under assumption A, the LSI fails, i.e. there is no $s>0$ such that
%\[
%s Ent_{\pi} (f) \leq {\cE}(\sqrt{f},\sqrt{f})
%\]
%for every $f \geq 0$.
%\end{proposition} 
%\begin{proposition}\label{pr3}
%In the Poisson case ($a(n) = \l$ and $b(n) = n$) the constant $c=1$ is the optimal constant in the MLSI.
%\end{proposition}
%Both Propositions \ref{pr2} and \ref{pr3} have a simple proof relying on test functions. For Proposition \ref{pr2} one takes the family of test functions $f_k(n) := {\bf 1}_{(k,+\infty)}(n)$, while for Proposition \ref{pr3} one chooses $f_k(n) := e^{-n/k}$. Details are simple and omitted.

The proof of Proposition \ref{nonmon} now follows by an application of
the perturbation argument recalled in Section \ref{rws}.
Namely, due to Lemma \ref{equiv} and Theorem \ref{t1} we know that
$\tilde \pi$ satisfies the inequality (\ref{entdisson}). Therefore
(\ref{entdisson}) follows from the equivalence between $\pi$ and $\tilde \pi$.

\subsection{Extension to $\bbZ$}
It is not difficult to extend the result of Theorem \ref{t1} to processes 
on $\bbZ$ rather than $\bbN$. Namely, consider the process with generator (\ref{genbd}) for all $n\in\bbZ$. Here, of course, we do not require $b(0)=0$. Again, we assume reversibility in the form (\ref{rever}), which holds now for every $n\in\bbZ$.
Similarly, we choose the function $R$ as in (\ref{rdb}) for all $n\in\bbZ$. It is easily checked that all the arguments given in the proof of Theorem \ref{t1} apply to this case without modification, provided the requirements of assumption ({\bf A}) are extended to all $n\in\bbZ$.  
For instance, this can be used to show that the double sided Poisson measures $\tilde\pi_\l(n) = (2\nep{\l}-1)^{-1}\l^{|n|}/|n|!$, $n\in\bbZ$, $\l\in(0,1)$ satisfy the inequality
$$
\ent_{\tilde \pi_\l}(f)\leq \frac{1}{1-\l}\,\cE(f,\log f)\,.
%\tilde\pi_\l\left[\nabla_+ f    \nabla_+ \log f \right]\,.
$$
Indeed, here we may choose $a(n)=\l$ for $n\geq 0$ and $b(n)=\l$ for $n\leq 0$. This gives $b(n)=n$ for all $n\geq 1$ and $a(n)=-n$ for all $n\leq -1$. In particular, $-\nabla_+ a(n) + \nabla_+b(n) = 1$ for all
$n\neq 0,-1$, in which cases it is equal to $c=1-\l$ so that Theorem \ref{t1} implies the above estimate. Several improvements of this type of estimates can be obtained along the lines discussed in the previous
subsections. On the other hand an extension to processes on
$\bbZ^d$, $d\geq 2$, does not appear to be straightforward.

\section{Zero range processes}

In this section we consider a class of interacting particle systems
consisting of finitely many particles moving in a finite set of
sites. Particles are neither created nor destroyed. The elements of
the set 
$\{1,2,\ldots,L\}$ label the sites; 
for $x \in \{1,2,\ldots,L\}$, $\eta_x \in \N$ 
denotes the number of particles at $x$. The whole configuration will
be denoted by $\eta \in S := \N^L$. 

The set $G$ of allowed moves is given by the set of maps from $S$ to $S$ of the form $\eta \mapsto \eta^{xy}$, with $x \neq y \in \{1,2,\ldots,L\}$, and
\[
\eta^{xy}_z = \left\{
\begin{array}{ll}
\eta_z & \mbox{if } z \not\in \{x,y\} \mbox{ or } \eta_x = 0  \\
\eta_x -1 & \mbox{for } z = x, \  \eta_x >0  \\
\eta_y +1 & \mbox{for } z = y ,  \  \eta_x >0.
\end{array}
\right.
\]
In other words $\eta^{xy}$ is obtained from $\eta$ by moving a particle (if any) from the site $x$ to the site $y$. We simply denote by $xy$ the map $\eta \mapsto \eta^{xy}$, and by $\nabla_{xy}$ the corresponding discrete gradient.

For $x \in \{1,2,\ldots,L\}$ consider  functions $c_x: \N \ra
[0,+\infty)$ such that $c_x(0) =0$, $c_x(n) >0$ for
$n>0$. $c_x(\eta_x)$ is the rate at which a particle is
moved from the site $x$ to a site $y$ chosen with uniform
probability. Thus we consider dynamics on $S$ for which the rate
$c(\eta,xy)$ of moving a particle from $x$ to $y$ is
$L^{-1}c(\eta_x)$. 
Therefore (\ref{geno}) becomes
\be{genzero}
\cL f(\eta) = \frac{1}{L} \sum_{x,y} c_x(\eta_x) \nabla_{xy} f(\eta)\,,
\end{equation}
where the sum extends to all $x,y\in\{1,\dots,L\}$. 
The continuous time Markov chain with 
generator (\ref{genzero}) is the zero--range process on the
complete graph with $L$ vertices. 

Note that the total number of particles $N:= \sum_x \eta_x$ is conserved. 
Set $p_x(n) := \prod_{k=1}^n \frac{1}{c_x(k)}$ for $n \geq 1$, $p_x(0)
=1$, and consider the probability $\pi_N$, defined on configurations
with $N$ particles, with $N =\sum_x \eta_x$, 
given by
\[
\pi_N(\eta) := \frac{1}{Z_N} \prod_{x=1}^L p_x(\eta_x),
\]
where $Z_N := \sum_{\eta \in S: \sum_x \eta_x = N}\prod_{x=1}^L
p_x(\eta_x)$ is the normalization. In what follows the  subscripts $N$
will be omitted. In the context of the class of models in Section
\ref{bbe}, we see easily that $(xy)^{-1} = yx$ and that the
reversibility condition ({\bf Rev}) holds, because of the identity
\be{revxy}
\pi\left[ c_x(\eta_x) g(\eta)
 \right] = \pi \left[c_y(\eta_y) g(\eta^{yx})\right]\,, 
\end{equation}
valid for arbitrary functions $g:S\to\bbR$.
We now define the function $R(\eta,\g,\d)$ to be
\be{rcll}
%L^2\,
R(\eta,xy,uv)  := \,\frac1{L^2}\,\begin{cases}
c_x(\eta_x) c_u(\eta_u) & \text{for } \;x\neq u\\
c_x(\eta_x) c_x(\eta_x -1)& \text{for } \;x = u
\end{cases}
\end{equation}
%
%
%
%\[
%\begin{array}{rcll}
%R(\eta,xy,uv) & := & \frac{1}{L^2}\,c_x(\eta_x) c_u(\eta_u) & \mbox{for } |\{x,y,u,v\}|=4 \\
%R(\eta,xy,xv) & := & \frac{1}{L^2}\,c_x(\eta_x) c_x(\eta_x -1)& \mbox{for } |\{x,y,v\}|=3 \\
%R(\eta,yx,xv) & := &  \frac{1}{L^2}\,c_y(\eta_y)c_x(\eta_x) & \mbox{for } |\{x,y,v\}|=3 \\
%R(\eta,xy,uy) & := &  \frac{1}{L^2}\,c_x(\eta_x) c_u(\eta_u) & \mbox{for } |\{x,y,v\}|=3 \\
%R(\eta,yx,uy) & := & \frac{1}{L^2}\,c_y(\eta_y) c_u(\eta_u) & \mbox{for } |\{x,y,v\}|=3 \\
%R(\eta,xy,xy) & := & \frac{1}{L^2}\,c_x(\eta_x) c_x(\eta_x -1)  & \\
%R(\eta,xy,yx) & := & \frac{1}{L^2}\,c_x(\eta_x) c_y(\eta_y) & 
%\end{array}
%\]
where $c(-1)$ is meant to be zero. The symmetry condition ({\bf P1}) is
checked immediately. Also condition ({\bf P3}) is simple to
check. Indeed, 
$xy$ and $uv$ commute when applied to $\eta$ unless $\eta_x
\eta_u =0$, but in this latter case $R(\eta,xy,uv)=0$. Condition ({\bf P2})
can be checked by direct inspection using (\ref{revxy}).
%amounts to check that $R(\eta,xy,uv) = R(\eta^{xy},yx,uv)$, which is
%checked by hands in 
%all cases \what. 

In order to use Proposition \ref{pr1} we shall assume:
\bi
\item[({\bf A})] All functions $c_x(\cdot)$ are nondecreasing.
\ei
\bl{ax}
Assume {\rm ({\bf A})}. Set $\G(\eta ,\g,\d) = c(\eta ,\g) c(\eta ,\d) - R(\eta
,\g,\d)$, where $c(\eta ,\g) = L^{-1} c_x(\eta_x)$, whenever $\g=xy$.
Then %have 
\begin{align}
 &\pi\left[ \sum_{\g,\d}\Gamma(x,\g,\d)\left( \nabla_{\g} f(\eta)
     \nabla_{\d} \log f(\eta) +  \frac{\nabla_{\g} f (\eta)
       \nabla_{\d} f(\eta)}{f(\eta)}\right) \right]\nonumber \\
&\quad\quad\quad\quad\quad\quad\geq \frac{1}{L}
 \sum_{x,y}\pi \left[ c_x(\eta_x) A_{x}(\eta) \nabla_{xy} f(\eta)
   \nabla_{xy} \log f(\eta) 
\right],\label{zr6}
 \end{align}
 where
 \[
 A_{x}(\eta) := \left(c_x(\eta_x ) - c_x(\eta_x-1)\right)\Big(1- \frac{1}{2L} \Big)
-\frac{1}{2L} \sum_{v:\;v\neq x} \left(c_v(\eta_v +1) - c_v(\eta_v) \right).
 \]
\el
\proof
We write
\begin{align}
&\pi\left[ \sum_{\g,\d}\Gamma(\eta,\g,\d)\left( \nabla_{\g} f(\eta)
    \nabla_{\d} \log f(\eta) +  \frac{\nabla_{\g} f (\eta) \nabla_{\d}
      f(\eta)}{f(\eta)}\right) \right] \nonumber\\ &\qquad 
%= \frac{1}{L^2}
%\sum_{xy,xv} \pi \left[ c_x(\eta_x) \left(c_x(\eta_x) -
%    c_x(\eta_x -1) \right) \left( \nabla_{xy} f(\eta) \nabla_{xv} \log
%    f(\eta) + \frac{ \nabla_{xy} f(\eta)  \nabla_{xv} f(\eta)
%    }{f(\eta)} \right) \right] \\
%\\ + \frac{1}{L^2} \sum_{xy} \pi \left[
%  c_x(\eta_x) \left(c_x(\eta_x) - c_x(\eta_x -1) \right) \left(
%    \nabla_{xy} f(\eta) \nabla_{xy} \log f(\eta) + \frac{ \nabla_{xy}
%      f(\eta)  \nabla_{xy} f(\eta) }{f(\eta)} \right) 
%\right] 
%\end{multline*}
%\begin{multline}
\label{zr1}
= \frac{1}{L^2} \sum_{x,y,v} \pi \left[ c_x(\eta_x) \left(c_x(\eta_x)
    - c_x(\eta_x -1) \right) \left( \nabla_{xy} f(\eta) \nabla_{xv}
    \log f(\eta) + \frac{ \nabla_{xy} f(\eta)  \nabla_{xv} f(\eta)
    }{f(\eta)} \right) \right] \nonumber\\ &\qquad 
\geq 
\frac{1}{L^2} \sum_{x,y,v} \pi \left[ c_x(\eta_x) \left(c_x(\eta_x) -
    c_x(\eta_x -1) \right) \nabla_{xy} f(\eta) \nabla_{xv} \log
  f(\eta) 
\right],
\end{align}
where in (\ref{zr1}) we used ({\bf A}) and the fact that, for every
$x$ and $\eta$:
\[
 \sum_{y,v} \nabla_{xy} f(\eta)  \nabla_{xv} f(\eta)  = \left[\sum_y \nabla_{xy} f(\eta) \right]^2 \geq 0\,.
 \]
 Now observe that
 \begin{align} \label{zr2}
 \frac{1}{L^2} &\sum_{x,y,v} \pi \left[ c_x(\eta_x) \left(c_x(\eta_x) -
     c_x(\eta_x -1) \right) \nabla_{xy} f(\eta) \nabla_{xv} \log
   f(\eta) \right] \nonumber\\ &\qquad
 =  \frac{1}{L^2} \sum_{x,y,v} \pi \left[
   c_x(\eta_x) \left(c_x(\eta_x) - c_x(\eta_x -1) \right) \nabla_{xy}
   f(\eta) \nabla_{xy} \log f(\eta) \right] \nonumber\\ &\qquad\quad
 + \frac{1}{L^2}
 \sum_{x,y,v} \pi \left[ c_x(\eta_x) \left(c_x(\eta_x) - c_x(\eta_x
     -1) \right) \nabla_{xy} f(\eta)\left( \log f(\eta^{xv}) - \log
     f(\eta^{xy}) \right) \right] \nonumber\\ &\qquad 
= \frac{1}{L} \sum_{x,y} \pi
 \left[ c_x(\eta_x) \left(c_x(\eta_x) - c_x(\eta_x -1) \right)
   \nabla_{xy} f(\eta) \nabla_{xy} \log f(\eta) \right] \nonumber\\ &\qquad\quad +
 \frac{1}{L^2} \sum_{x,y,v} \pi \left[ c_x(\eta_x) \left(c_x(\eta_x) -
     c_x(\eta_x -1) \right) f(\eta^{xy})\left( \log f(\eta^{xv}) -
     \log f(\eta^{xy}) \right) 
\right],
 \end{align}
 where in the last step we simply observed that, by symmetry, 
 \[
 \sum_{y,v} \left( \log f(\eta^{xv}) - \log f(\eta^{xy}) \right) = 0\,.
 \]
 We now use reversibility in the form (\ref{revxy})
to rewrite the last term in (\ref{zr2}):
 \begin{align}
   \frac{1}{L^2} &\sum_{x,y,v} \pi \left[ c_x(\eta_x)
     \left(c_x(\eta_x) - c_x(\eta_x -1) \right) f(\eta^{xy})\left(
       \log f(\eta^{xv}) - \log f(\eta^{xy}) \right) \right]
   \nonumber \\ 
& \quad = \frac{1}{L^2} \sum_{x,y,v:\;y\neq x} 
\pi \left[ c_y(\eta_y) \left(c_x(\eta_x +1) - c_x(\eta_x) \right)
  f(\eta)\left( \log f(\eta^{yv}) - \log f(\eta) \right) \right]\label{zr3} \\
& \quad \quad\quad+\frac{1}{L^2} \sum_{x,v} 
\pi \left[ c_x(\eta_x) \left(c_x(\eta_x ) - c_x(\eta_x-1) \right)
  f(\eta)\left( \log f(\eta^{xv}) - \log f(\eta) \right) \right]
\nonumber \\ 
& \quad =   \frac{1}{L^2} \sum_{x,y,v:\; v\neq x} 
\pi \left[ c_v(\eta_v) \left(c_x(\eta_x +1) - c_x(\eta_x) \right)
  f(\eta^{vy})\left( \log f(\eta) - \log f(\eta^{vy}) \right) \right]
\label{zr4}\\
& \quad\quad\quad + \frac{1}{L^2} \sum_{x,y} 
\pi \left[ c_x(\eta_x) \left(c_x(\eta_x ) - c_x(\eta_x-1) \right)
  f(\eta^{xy})\left( \log f(\eta) - \log f(\eta^{xy}) \right) \right]
\nonumber 
%
%
%\frac{1}{L^2} \sum_{x,y,v} \pi \left[
%  c_v(\eta_v) \left(c_x(\eta_x +1) - c_x(\eta_x) \right)
%  f(\eta^{vy})\left( \log f(\eta) - \log f(\eta^{vy}) \right) \right]
%\nonumber \\ & = & \frac{1}{L^2} \sum_{x,y,v} \pi \left[ c_y(\eta_y)
%  \left(c_x(\eta_x +1) - c_x(\eta_x) \right) f(\eta^{yv})\left( \log
%    f(\eta) - \log f(\eta^{yv}) \right) \right] \label{zr4} \\ & = & -
%\frac{1}{2L^2} \sum_{x,y,v}  \pi \left[ c_x(\eta_x) \left(c_v(\eta_v
%    +1) - c_v(\eta_v) \right) \nabla_{xy} f(\eta) \nabla_{xy} \log
%  f(\eta) 
%\right] \label{zr5}
 \end{align}
Therefore, exchanging the labels $y$ and $v$ in (\ref{zr4}) and
summing this expression with (\ref{zr3}) we obtain
\begin{align*}
  &  \frac{1}{L^2} \sum_{x,y,v} \pi \left[ c_x(\eta_x) \left(c_x(\eta_x) - c_x(\eta_x -1) \right) f(\eta^{xy})\left( \log f(\eta^{xv}) - \log f(\eta^{xy}) \right) \right]  \nonumber
\\
&\quad =  
- \frac{1}{2L^2} \sum_{x,y,v:\;y\neq x} 
\pi \left[ c_y(\eta_y) \left(c_x(\eta_x +1) - c_x(\eta_x) \right)
 \nabla_{yv} f(\eta)\nabla_{yv}\log f(\eta) \right]\label{zr5} \\
&\quad\quad\quad-\frac{1}{2L^2} \sum_{x,v} 
\pi \left[ c_x(\eta_x) \left(c_x(\eta_x ) - c_x(\eta_x-1) \right)
 \nabla_{xv} f(\eta)\nabla_{xv}\log f \right]\,.\nonumber
\end{align*}
%
% where (\ref{zr5}) is obtained by taking the arithmetic mean of (\ref{zr3}) and (\ref{zr4}) and then renaming the indexes. Inserting (\ref{zr5}) in (\ref{zr2}) and then in (\ref{zr1}) we obtain
% \be{zr6}
% \pi\left[ \sum_{\g,\d}\Gamma(x,\g,\d)\left( \nabla_{\g} f(\eta) \nabla_{\d} \log f(\eta) +  \frac{\nabla_{\g} f (\eta) \nabla_{\d} f(\eta)}{f(\eta)}\right) \right] \geq \frac{1}{L} \sum_{x,y}\pi \left[ c_x(\eta_x) A_{xy}(\eta) \nabla_{xy} f(\eta) \nabla_{xy} \log f(\eta) \right],
% \end{equation}
% where
% \[
% A_{xy}(\eta) := c_x(\eta_x +1) - c_x(\eta_x) - \frac{1}{2L} \sum_v \left(c_v(\eta_v +1) - c_v(\eta_v) \right).
% \]
The desired conclusion now follows from (\ref{zr1}) and (\ref{zr2}). \qed

The previous lemma allows us
 %This expression allows 
 to obtain ({\bf MLSI}) under the following condition:
 \bi
 \item[({\bf B})] There exist $0 \leq \d < c$ such that for every $x \in \{1,2,\ldots,L\}$ and $n \geq 0$
 \[
 c \leq c_x(n+1) - c_x(n) \leq c + \d.
 \]
 \ei
 Indeed, it is immediately seen that, under ({\bf B}), 
 \[
 A_{x}(\eta) \geq \frac{c-\d}{2}\,.
 \]
 Therefore, using Proposition \ref{pr1}, we have proved the following result. 
\bt{t2}
If the rates $c_x$ satisfy assumption {\rm ({\bf B})} then the inequality
(\ref{MLSI''}) holds with $\k=c-\d$, uniformly in the number of vertices and the number of particles. In particular, 
%for the corresponding zero--range process on the complete graph,
{\rm ({\bf MLSI})} holds with the same constant.
%$\a = c-\d$, uniformly in the number of vertices and the number of particles.
\et 

The following remarks give some elements to test the strength of 
the theorem we have just derived.

\subsection{The independent case}
Consider the case of linear rates, i.e.\ $c_x(\eta_x) = a_x
\eta_x$, for some coefficients $a_x\in(0,\infty)$,
$x\in\{1,\dots,L\}$. %Clearly, this case satisfies assumption {\bf (B)}
%provided that $c\leq a_x \leq c+\d$ for some $0<\d<c$. 
In this special case
%We observe that in the case of linear rates, i.e.\ $c_x(\eta_x) = a_x
%\eta_x$, for some coefficients $a_x\in(0,\infty)$,
%$x\in\{1,\dots,L\}$, 
the process describes $N$ independent random
walks on the complete graph, where each particle jumps from a vertex
$x$ to a vertex $y$ with rate $a_x$. 
The equilibrium measure $\pi$ becomes a
product of identical single--particle measures $\pi_1$, the
$\pi_1$--probability that the particle is at vertex $x$ being
proportional to $a_x^{-1}$. 
By the tensorization property of entropy one can then reduce the
problem to establishing ({\bf MLSI}) for a single random walk.
Already in this case, Theorem \ref{t2} gives a non--trivial result.
Note that, in the homogeneous case $a_x\equiv 1$, our estimate reduces
to the well known bound $\a=1$ for the simple random walk on the
unweighted complete graph, see e.g.\ Example 3.10 in \cite{BT}.

\subsection{Non--convex decay of entropy}
It is natural to wonder about the necessity of the restriction $\d<c$
in our assumption ({\bf B}). It was shown in \cite{bcdp} that as far
as the spectral gap is concerned, inequality ({\bf PI}) holds for this
model with $\g\geq c$ as soon as $c_x(\eta_x)-c_x(\eta_x-1)\geq c$ for
all $x$ and $\eta$, without further restriction. While we suspect that 
a similar condition should be sufficient for ({\bf MLSI}) it is
interesting to note that in order to have convexity of the relative entropy 
along the semigroup some restriction on the growth of the rates is
necessary. To see this we consider the following simple example of 
zero--range process exhibiting 
non--convex decay of entropy, i.e.\ such that
\be{nc}
\pi[\cL f \cL\log f] + \pi\left[ \frac{(\cL f)^2}{f} \right] < 0\,,
\end{equation}
for some $f>0$. 
Take $N=1$ particle only. Note that, since $N=1$ we must
have $R(\eta,xy,uv)=0$ for all $\eta$ and all $xy,uv$ in (\ref{rcll}).
Moreover, set $c_x:=c_x(\eta_x)$, $\pi_x:=\pi(\eta)$ and 
$f_x:=f(\eta)$ when the particle is at $x$ (i.e.\ when $\eta_x=1$).
Since $\pi_x=Z^{-1}c_x^{-1}$,  $Z:=\sum_x c_x^{-1}$, we 
see that the left hand side of (\ref{nc}) 
equals
\begin{align}
&\frac{1}{L^2}
\sum_{x,y,z}
\pi\left[c_{x}(\eta_x)^2 \left\{\nabla_{xy}f(\eta)\nabla_{xz}\log f(\eta)
+ \frac{\nabla_{xy}f(\eta)\nabla_{xz}f(\eta)}{f(\eta)}\right\}\right]
\nonumber\\
&\quad\quad = \frac{1}{Z\,L^2}
\sum_{x}
c_{x} \sum_{y,z}\left\{(f_y-f_x)\log(f_z/f_x)
+ \frac{(f_y-f_x)(f_z-f_x)}{f_x}\right\}\,%\quad\, q:=\big(\sum_x
                                %c_x^{-1}\big)^{-1}\,.
%\quad\quad 
=: \frac{1}{Z\,L^2}\sum_x c_x \,Q_x
\label{count1}
\end{align}
This expression can be shown to be negative for suitable choices of $\{f_x\}$
and $\{c_x\}$. A simple example is obtained if e.g.\ $L=3$, $f_1=1$, $f_2=2$,
$f_3=\e>0$ and $c_1>c_2=c_3=1$. If $\e$ is sufficiently small, in this
case we see
that $Q_1 = \e(\log 2+\e +\log \e) < 0$, so that 
$\sum_x c_x \,Q_x = c_1 Q_1 + Q_2 + Q_3$ must become negative when
$c_1$ is large. Thus (\ref{nc}) holds
and the entropy of $T_tf$ is not convex in $t\geq 0$.
Clearly, (\ref{count1}) can be used to construct many other examples
of such a behavior.
%On the other hand it is clear that such a system has
%a finite % (possibly depending on $L$) 
%entropy dissipation constant.
%In particular, 
%the entropy of $f_t=\exp{(t\cL)}f$ decays not slower than exponentially 
%but
%it does not so in a convex way. This is in striking contrast 
%with the case of
%the variance, where spectral theory shows that the decay is always
%convex.

\subsection{Non--monotone versus non--homogeneous rates}
The case of non--monotone rates refers to the situation where the
rates $c_x$ satisfy the assumptions appearing in Proposition
\ref{nonmon}. Unfortunately, Theorem \ref{t2} does not extend 
to this case by simple perturbation arguments.
Zero range processes with non--monotone rates have 
been thoroughly studied in the literature,
under
the further assumption that the model is homogeneous, i.e.\ $c_x=c_y$
for all $x,y$. 
For the nearest neighbor version of this model, both
Poincar\'e and logarithmic Sobolev inequalities have been established
\cite{LSV,DP}. Moreover, the corresponding complete graph model
has been shown to satisfy the ({\bf MLSI}) inequality \cite{CP}. 
These results, all based on some version of the so--called 
martingale decomposition method, do not extend 
to non--homogeneous models in a standard way and, 
as far as we know Theorem \ref{t2}
represents the only criterium available in
non--homogeneous models.

\section{Bernoulli-Laplace models}
 
 As in previous section we consider a system of particles moving in
 the finite set of sites $\{1,2,\ldots,L\}$; here we assume that particles
 are subject to an exclusion rule, namely at most one particle per
 site is allowed. Thus $S := \{0,1\}^L$. The set of allowed moves
 is $G := \{ xy:\;x,y \in \{1,2,\ldots,L\}\,,\; x \neq y\}$, where, for
 $\eta \in S$, $\eta^{xy} = \eta$ unless 
$\eta_x (1-\eta_y) =1$, and in this case
 \[
\eta^{xy}_z = \left\{
\begin{array}{ll}
\eta_z & \mbox{if } z \not\in \{x,y\}   \\
0 & \mbox{for } z = x,  \\
1 & \mbox{for } z = y .
\end{array}
\right.
\]
To each site $x$ we associate a Poisson clock of constant 
intensity $\l_x>0$; when the
clock of site $x$ rings, a site $y$ is chosen at random: if
$\eta_x =1$ and $\eta_y =0$ then the particle at $x$ moves to $y$,
otherwise nothing happens. This dynamics corresponds 
to the infinitesimal generator
\[
\cL f(\eta) := \frac{1}{L} \sum_{x,y =1}^L \l_x \eta_x (1-\eta_y) \nabla_{xy} f(\eta).
\]
In other words, we set $c(\eta,xy) = L^{-1}\l_x \eta_x (1-\eta_y)$ in
(\ref{geno}). Denote by
$N \leq L$ the number of particles 
in the system; since it is conserved by the dynamics, we can consider
the restriction of the dynamics to configurations with $N$
particles. In this restricted state space there is a unique stationary
distribution $\pi_N$, given by conditioning to configurations with $N$ particles the product of Bernoulli measures with parameters $\frac{1}{1+\l_x}$. More precisely
\[
\pi_N(\eta) = \frac{1}{Z_{L,N} }\prod_{x=1}^L \left(\frac{1}{1+\l_x}\right) ^{\eta_x} \left(\frac{\l_x}{1+\l_x}\right) ^{1-\eta_x} .
\]
We will refer to this measure as the {\em canonical} measure on $\{1,2,\ldots,L\}$.
The reversibility condition ({\bf Rev}) holds true for $\pi_N$; the subscript $N$ will be omitted from now on.
Note that, for $N=1$, zero range processes coincide with Bernoulli-Laplace models. In particular, the counterexample in Section 4.2 concerning the non convex decay of entropy applies here too. Thus, some bound on the non homogeneity of the model is needed. The following condition, which is presumably not optimal, is analogous to condition ({\bf B}) for zero range processes:
\bi
\item[({\bf B})] 
There exists $0\leq \d < c$ such that for every $x=1,2,\ldots,L$
\[
c \leq \l_x \leq c+\d.
\]
\ei

\bt{thBL} 
Assume {\rm ({\bf B})}. Then the inequality (\ref{MLSI''}) holds with $\k = c-\d$. In particular, 
{\rm ({\bf MLSI})} holds with $\a  = c-\d$.
\et

It should be stressed that because in the present case only one particle per site is allowed, the
proof of the next Theorem~\ref{thBL} requires some arguments that were not needed the proof of Theorem~\ref{t2}.

\smallskip\noindent

{\bf Proof of Theorem \ref{thBL}}. 
We use, of course, Proposition \ref{pr1}. The choice of $R$ is essentially forced by the commutation condition ({\bf P3}). We set
\[
R(\eta,xy,zu) := \frac1{L^2}\left\{ \begin{array}{ll} \l_x \l_z \eta_x (1-\eta_y) \eta_z (1-\eta_u) & \mbox{for } |\{x,y,z,u\}| =4 \\ 0 & \mbox{otherwise.} \end{array} \right.
\]
We obtain, after having noticed that $\nabla_{xy} f \nabla_{yz} g
\equiv 0 \equiv \nabla_{xy}f \nabla_{ux} g$ for any choice of
$x,y,z,u$ and $f,g$,
\begin{align}\label{BL1} 
&\pi\left[ \sum_{\g,\d}\Gamma(\eta,\g,\d)\left( \nabla_{\g} f(\eta)
    \nabla_{\d} \log f(\eta) +  \frac{\nabla_{\g} f (\eta) \nabla_{\d}
      f(\eta)}{f(\eta)}\right) \right] 
\nonumber\\ &\quad
= 
\frac{1}{L^2} \sum_{\stackrel{ x,y,z:}{\scriptscriptstyle
    |\{x,y,z\}|=3}} \pi \left[\l_x^2 \nabla_{xy}f \nabla_{xz} \log f 
\right] 
+ \frac{1}{L^2} \sum_{\stackrel{ x,y,u:}{\scriptscriptstyle
    |\{x,y,u\}|=3}} \pi \left[\l_x \l_u \nabla_{xy}f \nabla_{uy} \log
  f \right]  
\nonumber\\ &\qquad\quad
+\frac{1}{L^2} \sum_{x,y}  \pi \left[\l_x^2 \nabla_{xy}f
  \nabla_{xy} \log f \right] 
+ \frac{1}{L^2} \sum_{\stackrel{ x,y,u:}{\scriptscriptstyle
    |\{x,y,u\}|=3}} \pi \left[\l_x \l_u \frac{\nabla_{xy}f \nabla_{uy}
    f }{f}
\right]  
\nonumber\\ &\qquad\quad
+ \frac{1}{L^2} \sum_{\stackrel{ x,y,z:}{\scriptscriptstyle
    |\{x,y,z\}|=3}} \pi \left[\l_x^2 \frac{\nabla_{xy}f \nabla_{xz} f
  }{f}
\right] 
+ \frac{1}{L^2} \sum_{x,y}  \pi \left[\l_x^2 
\frac{\left(\nabla_{xy}f
    \right)^2}{f} 
\right].
\end{align}
The sum of the last two terms in (\ref{BL1}) equals
\be{BL2}
\frac{1}{L^2}\sum_x \l_x^2  \pi\left[ \frac{\left(\sum_y \nabla_{xy} f \right)^2}{f} \right] \geq 0 \,.
\end{equation}
We now claim that
\be{BL3}
\sum_{\stackrel{ x,y,u:}{\scriptscriptstyle |\{x,y,u\}|=3}} \pi \left[\l_x \l_u \frac{\nabla_{xy}f \nabla_{uy}  f }{f}\right]   \geq 0\,.
\end{equation}
for every $f>0$. To prove (\ref{BL3}) we observe that
\begin{align*}
\sum_{\stackrel{ x,y,u:}{\scriptscriptstyle |\{x,y,u\}|=3}} \pi
\left[\l_x \l_u \frac{\nabla_{xy}f \nabla_{uy}  f }{f}\right] &= \sum_T
\sum_{\stackrel{ x,y,z \in T:}{\scriptscriptstyle |\{x,y,z\}|=3}} \pi
\left[\l_x \l_u \frac{\nabla_{xy}f \nabla_{uy}  f }{f}\right] \\& =
\sum_T \sum_{\stackrel{ x,y,z \in T:}{\scriptscriptstyle
    |\{x,y,z\}|=3}} \pi \left\{ \pi\left[\l_x \l_u \frac{\nabla_{xy}f
      \nabla_{uy}  f }{f}\bigg|\eta_{T^c}\right] 
\right\},
\end{align*}
where %in $\sum_T$ 
$T$ varies over $\{T \subseteq \{1,2,\ldots,L\}: |T|=3\}$,
and $\pi[\, \cdot \, \tc\eta_{T^c}]$ denotes the conditional
expectation with respect to the configuration outside $T$, which we denote by $\eta_{T^c}$. Note that
the corresponding conditional measure is the canonical measure on $T$  with $N - \sum_{x \not\in T} \eta_x$
particles. Inequality (\ref{BL3}) is proved if we show that for every
fixed $T$ and 
$\eta_{T^c}$
\be{BL4}
\sum_{\stackrel{ x,y,z \in T:}{\scriptscriptstyle |\{x,y,z\}|=3}}\pi\left[\l_x \l_u \frac{\nabla_{xy}f \nabla_{uy}  f }{f}\bigg|\eta_{T^c}\right] \geq 0.
\end{equation}
Set $T = \{a,b,c\}$ note that if $\sum_{x \not\in T} \eta_x \in \{N,
N-3\}$ than $T$ contains either no particle or no hole, so that
$\nabla_{xy}f \equiv 0$ for every $f$ and every $x,y \in
T$. Similarly, for $\sum_{x \not\in T} \eta_x = N-1$ there is only one
particle in $T$, so $\nabla_{xy}f \nabla_{uy}  f \equiv 0$ for $x \neq
u$. So we only need to consider the case $\sum_{x \not\in T} \eta_x =
N-2$, which means there is exactly one hole in $T$. For a given 
configuration $\eta_{T^c}$ outside of $T$, we denote by
$\a$ the value of $f$ on the configuration with the hole in $a$. 
Similarly, $\b$ is the value of $f$
when the hole is in $b$ and $\g$ when the hole is in $c$. 
Then, by direct computation we get
\begin{align}
\frac{Z_{T,2}}{2}&\sum_{\stackrel{ x,y,z \in T:}{\scriptscriptstyle
    |\{x,y,z\}|=3}}
\pi\left[\l_x \l_u \frac{\nabla_{xy}f \nabla_{uy}  f }{f}\bigg|
  \eta_{T^c}\right] = \frac{\l_a}{1+\l_a} \frac{1}{1+\l_b} \frac{1}{1+\l_c}\l_b \l_c \frac{(\b-\a)(\g - \a)}{\a} \nonumber \\ & ~~~~~~ + \frac{\l_b}{1+\l_b} \frac{1}{1+\l_a} \frac{1}{1+\l_c}\l_a \l_c \frac{(\g-\b)(\a -
  \b)}{\b} %\nonumber \\ & ~~~~~~~
+ \frac{\l_c}{1+\l_c} \frac{1}{1+\l_a} \frac{1}{1+\l_b}\l_a \l_b  \frac{(\b-\g)(\a - \g)}{\g} \nonumber \\ 
& \quad\quad\qquad
= \frac{\l_a \l_b \l_c}{(1+\l_a)(1+\l_b)(1+\l_c)} \left[\frac{\b \g}{\a} +
\frac{\a \g}{\b} + \frac{\a \b}{\g} -\a - \b - \g \right] \,.
 \label{BL5}%\quad\quad\quad\quad
\end{align}
 where $Z_{T,2}$ is the normalization factor of the canonical measure on $T$ with $2$ particles, and the further factor $1/2$
 in the l.h.s.\ of (\ref{BL5}) is due to the fact that the sum over the
 ``particles'' $x,u$ is a sum over ordered pairs. We need to show that the
 r.h.s.\ of (\ref{BL5}) is nonnegative, for every $\a,\b,\g>0$. Since
 the expression is homogeneous of degree one and invariant for
 permutations of the variables, we may restrict to $\g=1$, $\a,\b \geq
 1$. In other words we need to 
show that
 \be{BL6}
 F(\a,\b) := \frac{\b }{\a} + \frac{\a }{\b} + \a \b -\a - \b - 1 \geq 0\,,
 \end{equation}
 for every $\a,\b \geq 1$. Since $z+z^{-1}\geq 2$ for all $z>0$, we
 have 
 \[
 F(\a,\b)\geq 1+\a \b -\a - \b = (\a-1)(\b-1) \geq 0,
 \]
 and (\ref{BL6}), follows. This shows (\ref{BL4}).
% 
% $F(\a,\b)\geq 1+\a \b -\a - \b$. Setting $\a=1+u,\b=1+v$, we
% see that $F(\a,\b)\geq uv$, which implies the claim.
%This proves (\ref{BL4}). 
%
%
% Note that $F(1,1) = 0$. To prove that $(1,1)$ is an absolute minimum, it is enough to show that, for every $m \geq 0$
% \be{BL6}
% \frac{d}{d\a} F(\a, 1+m(\a-1)) \geq 0
% \end{equation}
% for every $\a \geq 1$, i.e. $F$ is non-increasing along the lines through $(1,1)$ with positive slope, and that
% \be{BL7}
% \frac{d}{d\b} F(1,\b) \geq 0
% \end{equation}
% for $\b \geq 1$, to include the line with infinite slope. Let us show (\ref{BL6}).
% \begin{multline*}
%  \frac{d}{d\a} F(\a, 1+m(\a-1)) =  \frac{d}{d\a} \left[ \frac{1+m(\a-1) }{\a} + \frac{\a }{1+m(\a-1)} + \a [1+m(\a-1)] -\a - [1+m(\a-1)] \right. \\ (m-1) \left[\frac{1}{\a^2} - \frac{1}{[1+m(\a-1) ]^2} \right] + 2m(\a-1).
%  \end{multline*}
%Both summands in this last expression are nonnegative, since $\a \geq
%1$ and $\frac{1}{\a^2} - \frac{1}{[1+m(\a-1) ]^2}$ has the same sign
%of $(m-1)$, as it is shown easily. This establishes
%(\ref{BL6}). Inequality (\ref{BL7}) is easier, and the proof is 
%omitted.

From (\ref{BL1})-(\ref{BL4}) we get
\begin{align} 
&\pi\left[ \sum_{\g,\d}\Gamma(x,\g,\d)\left( \nabla_{\g} f(\eta)
    \nabla_{\d} \log f(\eta) +  \frac{\nabla_{\g} 
f (\eta) \nabla_{\d} f(\eta)}{f(\eta)}\right) \right] \nonumber\\ 
&\quad\quad\geq 
\frac{1}{L^2} \sum_{\stackrel{ x,y,z:}{\scriptscriptstyle
    |\{x,y,z\}|=3}} \pi \left[\l_x^2 \nabla_{xy}f \nabla_{xz} 
\log f \right] 
+ \frac{1}{L^2} \sum_{\stackrel{ x,y,u:}{\scriptscriptstyle
    |\{x,y,u\}|=3}} \pi \left[\l_x \l_u \nabla_{xy}f \nabla_{uy} \log f \right]
\nonumber\\ & \qquad\qquad\qquad
+ \frac{1}{L^2} \sum_{x,y}  \pi \left[\l_x^2 \nabla_{xy}f \nabla_{xy} \log f
\right] \nonumber\\ 
&\quad\quad = \frac{1}{L^2} \sum_{ x,y,z} \pi \left[\l_x^2 \nabla_{xy}f
  \nabla_{xz} 
\log f \right]\ 
+ \frac{1}{L^2} \sum_{\stackrel{ x,y,u:}{\scriptscriptstyle
    |\{x,y,u\}|=3}} \pi \left[\l_x \l_u \nabla_{xy}f \nabla_{uy} 
\log f \right]  \,.\label{BL8}
\end{align}
We now deal separately with the last two terms in (\ref{BL8}), 
similarly to what done in Section 4. For the first term we have
\begin{align} \label{BL9}
\sum_{ x,y,z} \pi \left[\l_x^2 \nabla_{xy}f \nabla_{xz} \log f \right] &=
\sum_{ x,y,z} \pi \left[\l_x^2 \eta_x (1-\eta_y) (1-\eta_z)\nabla_{xy}f
  \nabla_{xz} \log f \right] \nonumber\\
& 
= \sum_{ x,y,z} \pi \left[\l_x^2 \eta_x
  (1-\eta_y) (1-\eta_z)\nabla_{xy}f \nabla_{xy} \log f \right] \nonumber\\
&\qquad
+ \sum_{ x,y,z} \pi \left[\l_x^2 \eta_x (1-\eta_y) (1-\eta_z)\nabla_{xy}f(\eta)
  [\log f(\eta^{xz}) - \log f(\eta^{xy})] \right] \nonumber\\ & 
= (L-N)
\sum_{x,y} \pi \left[\l_x^2 \nabla_{xy}f \nabla_{xy} \log f \right]
\nonumber
\\ &\qquad+ \sum_{
  x,y,z} \pi \left[\l_x^2 \eta_x (1-\eta_y) (1-\eta_z)f(\eta^{xy}) [\log
  f(\eta^{xz}) - \log f(\eta^{xy})] \right] \,,
\end{align}
where we use the fact that, for any $z$ and $\eta$, 
$$
\sum_{z} (1-\eta_z) = L-N\,,
$$ 
and the fact that, for any $x$ and $\eta$: 
$$
\sum_{ y,z} (1-\eta_y) (1-\eta_z)
  [\log f(\eta^{xz}) - \log f(\eta^{xy})] = 0\,.
$$
Using reversibility we have
\begin{align}
&\sum_{ x,y,z} \pi\left[\l_x^2 \eta_x (1-\eta_y) (1-\eta_z)f(\eta^{xy}) [\log
  f(\eta^{xz}) - \log f(\eta^{xy})] \right]  \nonumber\\ 
&\quad
=  \sum_{ x,y,z:\,z\neq x} \pi \left[\l_x \l_y \eta_y (1-\eta_x)
  (1-\eta_z)f(\eta) \nabla_{yz}\log f(\eta) \right] \nonumber \\ 
& \quad = 
- \sum_{ x,y,z:\,y\neq x} \pi \left[\l_x \l_z \eta_z (1-\eta_y)
  (1-\eta_x)f(\eta^{zy}) \nabla_{zy}\log f(\eta) \right]\nonumber \\ 
& \quad = 
-\frac{1}{2}\sum_{ x,y,z:\,z\neq x} \pi \left[\l_x \l_y \eta_y (1-\eta_x)
  (1-\eta_z)\nabla_{yz}f \nabla_{yz} \log f \right]\nonumber\\
%& \quad
% = -\frac{(L-N-1)}{2}
%\sum_{y,z} \pi \left[\eta_y 
%(1-\eta_z)\nabla_{yz}f \nabla_{yz} \log f \right]\,,
\label{BL10}
\end{align}
where we use permutation of $x,y,z$. By condition ({\bf B}) we easily have
\[
\sum_{ x:\,z\neq x}\l_x  (1-\eta_x)(1-\eta_z) \leq (c+\d) (L-N-1)(1-\eta_z) \,,
\]
which, together with (\ref{BL9}) and (\ref{BL10}) yields
\be{BL10.5}
\sum_{ x,y,z} \pi \left[\l_x^2 \nabla_{xy}f \nabla_{xz} \log f \right] \geq \left[c(L-N) - \frac{c+\d}{2}(L-N-1) \right]\sum_{x,y} \pi \left[\l_x \nabla_{xy}f \nabla_{xy} \log f \right]
\,.
\end{equation}
The last summand in (\ref{BL8}) is dealt with similarly:
\begin{align}
\label{BL11}
\sum_{\stackrel{ x,y,u:}{\scriptscriptstyle |\{x,y,u\}|=3}} &\pi
\left[\l_x \l_u \nabla_{xy}f \nabla_{uy} \log f \right] = \sum_{\stackrel{
    x,y,u:}{\scriptscriptstyle |\{x,y,u\}|=3}} \pi \left[\l_x \l_u \eta_x \eta_u
  (1-\eta_y)\nabla_{xy}f \nabla_{uy} \log f \right] \nonumber \\ &=
\sum_{\stackrel{ x,y,u:}{\scriptscriptstyle |\{x,y,u\}|=3}} \pi
\left[\l_x \l_u \eta_x \eta_u (1-\eta_y)\nabla_{xy}f \nabla_{xy} \log f \right]
\nonumber \\ & 
\qquad\qquad+ \sum_{\stackrel{ x,y,u:}{\scriptscriptstyle |\{x,y,u\}|=3}} \pi
\left[\l_x \l_u \eta_x \eta_u (1-\eta_y)f(\eta^{xy})[\log f(\eta^{uy}) - \log
  f(\eta^{xy})]\right] \nonumber \\ &\geq  c(N-1) \sum_{x,y} \pi \left[\l_x \nabla_{xy}f
  \nabla_{xy} \log f \right] + \sum_{\stackrel{
    x,y,u:}{\scriptscriptstyle |\{x,y,u\}|=3}} \pi \left[ \l_y \l_u \eta_y
  (1-\eta_x)\eta_u f(\eta) \nabla_{ux} \log f(\eta) \right] \nonumber
\\ &
= c(N-1)
\sum_{x,y} \pi \left[\l_x \nabla_{xy}f \nabla_{xy} \log f \right] -
\sum_{\stackrel{ x,y,u:}{\scriptscriptstyle |\{x,y,u\}|=3}} \pi
\left[\l_y \l_x \eta_y \eta_x (1-\eta_u) f(\eta^{xu}) \nabla_{xu} \log f(\eta)
\right] \nonumber \\ &= c (N-1) \sum_{x,y} \pi \left[\l_x \nabla_{xy}f \nabla_{xy} \log f
\right] - \frac{1}{2} \sum_{\stackrel{ x,y,u:}{\scriptscriptstyle
    |\{x,y,u\}|=3}} \pi \left[ \l_x \l_u \eta_x \eta_u (1-\eta_y) \nabla_{xy}f
  \nabla_{xy} \log f \right] \nonumber \\& \geq \left[ c(N-1) - \frac{c+\d}{2} (N-1) \right]\sum_{x,y} \pi
\left[\l_x \nabla_{xy}f \nabla_{xy} 
\log f \right]\,.%\qquad\qquad\qquad\quad\quad\quad\quad\qquad
\end{align}
By (\ref{BL8}), (\ref{BL10.5}) and (\ref{BL11}) we end up with
\[
\pi\left[ \sum_{\g,\d}\Gamma(x,\g,\d)\left( \nabla_{\g} f(\eta)
    \nabla_{\d} \log f(\eta) +  \frac{\nabla_{\g} f (\eta) \nabla_{\d}
      f(\eta)}{f(\eta)}\right) \right] \geq \frac{c-\d}{2L}
\sum_{x,y} \pi \left[\l_x \nabla_{xy}f \nabla_{xy} 
\log f \right]
\]
which, together with Proposition \ref{pr1}, completes the proof.
\qed

\smallskip

We conclude with a comparison of our bound with previously known
results, that are limited to the homogeneous case $\l_x \equiv 1$. In this case we can take $c=1$ and $\d=0$ in Theorem \ref{thBL}, therefore obtaining the ({\bf MLSI}) with $\a=1$ . By using Yau's martingale method, Gao and Quastel \cite{GQ} and Goel
\cite{Goel} have
proved ({\bf MLSI}) for this model with $\a = 1/2$ (asymptotically as
$L\to\infty$). The same estimate
has been obtained with a similar approach by Bobkov and Tetali \cite{BT}. 
We mention that a proof of our bound $\a= 1$ in the homogeneous case 
can be also 
obtained by a slightly different method \cite{CT}.

It is not hard to show
that the homogeneous Bernoulli-Laplace model has a spectral gap equal to $1$, independent of
the number of particles (see e.g.\
\cite{bcdp}). Therefore, from (\ref{comp}) 
we have that the best constant $\a$ in ({\bf MLSI})
satisfies $\a\in[1,2]$. The optimal value is not known, even in the
case of one particle (random walk on the complete graph), see e.g.\
the discussion of Example
3.10 in \cite{BT}. Moreover,
it is not known whether the optimal constant depends on the number of
particles. 

For the non--homogeneous model considered here only the spectral gap
has been obtained before, see \cite{Cshonan}, where a uniform
Poincar\'e inequality is established under the assumption that
$C^{-1}\leq \l_x \leq C$ for some $C>0$, for all $x$.
 
%Moreover, standard
%comparison with spectral gap implies that the largest constant in MLSI
%is at most $2$. 
%We do not know whether the ``correction''
%$\frac{1}{L}$ in Theorem \ref{thBL} is substantially related to the
%fact that inequality (\ref{MLSI''}) is stronger that MLSI, 
%but we honestly doubt it. 

\end{document}